\documentclass[11pt]{article}
\usepackage{amsthm,amssymb,amsmath,amscd}
\usepackage[all]{xy}
\usepackage{mathrsfs}
\usepackage{multirow}

\usepackage{mathptm}

\newtheorem{Def}{Definition}[section]
\newtheorem{Prop}[Def]{Proposition}
\newtheorem{Theo}[Def]{Theorem}
\newtheorem{Lem}[Def]{Lemma}
\newtheorem{Koro}[Def]{Corollary}

\newcommand{\add}{{\rm add}}
\newcommand{\con}{{\rm con}}
\newcommand{\domdim}{{\rm dom.dim}}
\newcommand{\Hom}{{\rm Hom }}
\newcommand{\pd}{{\rm pd }}
\newcommand{\rad}{{\rm rad}}
\newcommand{\soc}{{\rm soc}}

\newcommand{\StHom}{{\rm \underline{Hom} }}
\newcommand{\repdim}{{\rm rep.dim }}
\newcommand{\gldim}{{\rm gl.dim}}
\newcommand{\findim}{{\rm fin.dim}}
\newcommand{\End}{{\rm End}}
\newcommand{\Ext}{{\rm Ext}}

\newcommand{\cpx}[1]{#1^{\bullet}}
\newcommand{\D}[1]{\mathscr{D}(#1)}
\newcommand{\Dz}[1]{{\mathscr{D}}^+(#1)}
\newcommand{\Df}[1]{{\mathscr{D}}^-(#1)}
\newcommand{\Db}[1]{ \mathscr{D}^{\rm b}(#1)}
\newcommand{\C}[1]{\mathscr{C}(#1)}
\newcommand{\Cz}[1]{{\mathscr{C}}^+(#1)}
\newcommand{\Cf}[1]{{\mathscr{C}}^-(#1)}
\newcommand{\Cb}[1]{{\mathscr{C}^b}(#1)}
\newcommand{\K}[1]{\mathscr{K}(#1)}
\newcommand{\Kz}[1]{\mathscr{K}^+(#1)}
\newcommand{\Kf}[1]{\mathscr{K}^-(#1)}
\newcommand{\Kb}[1]{ \mathscr{K}^{\rm b}(#1)}
\newcommand{\modcat}[1]{#1\mbox{{\rm -mod}}}
\newcommand{\stmodcat}[1]{#1\mbox{{\rm -{\underline{mod}}}}}
\newcommand{\pmodcat}[1]{#1\mbox{{\rm -proj}}}
\newcommand{\imodcat}[1]{#1\mbox{{\rm -inj}}}
\newcommand{\opp}{^{\rm op}}
\newcommand{\otimesL}{\otimes^{\rm\bf L}}
\newcommand{\otimesP}{\otimes^{\bullet}}
\newcommand{\rHom}{{\rm\bf R}{\rm Hom}}

\newcommand{\HomP}{{\rm Hom}^{\bullet}}

\newcommand{\lra}{\longrightarrow}
\newcommand{\ra}{\rightarrow}
\newcommand{\dickebox}{{\vrule height5pt width5pt depth0pt}}

\begin{document}
{\Large \bf
\begin{center}
Derived equivalences and stable equivalences of Morita type, I.
\end{center}}
\medskip

\centerline{\bf Wei Hu and Changchang Xi$^*$}
\begin{center} School of Mathematical Sciences, \\
Laboratory of Mathematics and Complex Systems, Beijing Normal University, \\
100875 Beijing, People's Republic of  China \\E-mail: hwxbest@163.com \quad  xicc@bnu.edu.cn\\
\end{center}

\medskip
\centerline{\emph{Dedicated to Professor Shaoxue Liu on the occasion
of his 80th Birthday}}

\begin{abstract}
For self-injective algebras, Rickard proved that each derived
equivalence induces a stable equivalence of Morita type. For general
algebras, it is unknown when a derived equivalence implies a stable
equivalence of Morita type. In this paper, we first show that each
derived equivalence $F$ between the derived categories of Artin
algebras $A$ and $B$ arises naturally a functor $\bar{F}$ between
their stable module categories, which can be used to compare certain
homological dimensions of $A$ with that of $B$; and then we give a
sufficient condition for the functor $\bar{F}$ to be an equivalence.
Moreover, if we work with finite-dimensional algebras over a field,
then the sufficient condition guarantees the existence of a stable
equivalence of Morita type. In this way, we extend the classic
result of Rickard. Furthermore, we provide several inductive methods
for constructing those derived equivalences that induce stable
equivalences of Morita type. It turns out that we may produce a lot
of (usually not self-injective) finite-dimensional algebras which
are both derived-equivalent and stably equivalent of Morita type,
thus they share many common invariants.
\end{abstract}

\renewcommand{\thefootnote}{\alph{footnote}}
\setcounter{footnote}{-1} \footnote{ $^*$ Corresponding author.
Email: xicc@bnu.edu.cn; Fax: 0086 10 58802136; Tel.: 0086 10
58808877.}
\renewcommand{\thefootnote}{\alph{footnote}}
\setcounter{footnote}{-1} 
\renewcommand{\thefootnote}{\alph{footnote}}
\setcounter{footnote}{-1} 

\section{Introduction}
As is well-known, derived equivalence and stable equivalence of
Morita type are two fundamental types of equivalences in algebras
and categories, and play an important role in the modern
representation theory of groups and algebras, they transfer
information from one algebra to another, and provide a convenient
bridge between two different (derived or stable) categories. In
particular, derived equivalences preserve many significant
invariants; for example, the center of an algebra, the number of
non-isomorphic simple modules, the Hochschild cohomology groups and
Cartan determinants, while stable equivalences of Morita type,
introduced in around 1990 (see \cite{Broue1994}, for example)  and
appearing frequently in the block theory of finite groups, preserve
also many nice invariants; for instance, the global, finitistic, and
representation dimensions \cite{Xi1} as well as the representation
types \cite{KrauseRpType}. For self-injective algebras, the two
notions are closely related to each other, this was revealed by a
well-known result of Rickard \cite{RickDFun}, which states that a
derived equivalence between self-injective algebras induces always a
stable equivalence of Morita type. Moreover, the remarkable Abelian
Defect Group Conjecture of Brou\'e, which states that the module
categories of a block algebra $A$ of a finite group algebra and its
Brauer correspondent $B$ should have equivalent derived categories
if their common defect group is abelian (see \cite{R3}), makes the
two concepts more attractive and intimate. However, for general
finite-dimensional algebras, derived equivalence and stable
equivalence of Morita type seem to be completely different from each
other; for example, a representation-finite algebra may be
derived-equivalent to a representation-infinite algebra via a
tilting module, and consequently they neither are stably equivalent
of Morita type nor have the same representation dimension. Thus a
natural question arises: What kind of relationship between a derived
equivalence and a stable equivalence of Morita type for general
finite-dimensional algebras could exist? In other words, we consider
the following question:

\smallskip
{\bf Question.} When does a derived equivalence between two
finite-dimensional (not necessarily self-injective) algebras $A$ and
$B$ induces a stable equivalence of Morita type between them ?

\smallskip
Thus, a positive answer to the above questions would let us know
more invariants between algebras $A$ and $B$. However, in the past
time, little is known about this question. One even does not known
when a derived equivalence induces a stable equivalence for general
finite-dimensional algebras.

In the present paper, we shall provide some answers to this
question. To state our main result, let us introduce the notion of
an almost $\nu$-stable functor. Suppose $F$ is a derived equivalence
between two Artin algebra $A$ and $B$, with the quasi-inverse
functor $G$. Further, suppose
$$\cpx{T}: \quad \cdots \lra 0\lra T^{-n}\lra
\cdots \lra T^{-1}\lra T^0\lra 0\lra \cdots$$ is a radical tilting
complex over $A$ associated to $F$, and suppose
$$\cpx{\bar{T}}: \quad \cdots \lra 0\lra \bar{T}^{0}
\lra \bar{T}^{1}\lra \cdots \lra \bar{T}^n\lra 0\lra \cdots
$$ is a radical tilting complex over $B$ associated to $G$. The functor $F$ is called {\it almost
$\nu$-stable} if $\add(\bigoplus_{i=-1}^{-n}
T^i)=\add(\bigoplus_{i=-1}^{-n} \nu_A T^i)$, and
$\add(\bigoplus_{i=1}^{n} \bar{T}^i)=\add(\bigoplus_{i=1}^n
\nu_B\bar{T}^i)$, where $\nu_A$ is the Nakayama functor of $A$. Note
that the summations exclude only the term in degree 0. If $A$ and
$B$ are self-injective, every derived equivalence between $A$ and
$B$ is almost $\nu$-stable (by Proposition \ref{prop4.1} below).
Surprisingly, even beyond the class of self-injective algebras there
are plenty of almost $\nu$-stable derived equivalences, for example,
the derived equivalences constructed in \cite[Corollary 3.8]{HuXi2}
and in Proposition \ref{stabletilting} below. In fact, we shall give
a general machinery below to produce such derived equivalences.

With this notion in mind, our main result can be stated as follows.

\begin{Theo}
Let $A$ and $B$ be Artin algebras, and let $F$ be a derived
equivalence between $A$ and $B$. Then:

$(1)$ $F$ induces a functor $\bar{F}$ from the stable module
category over $A$ to that over $B$.

$(2)$ If $F$ is almost $\nu$-stable, then the functor $\bar{F}$
defined in $(1)$ is an equivalence. Furthermore, if $F$ is an almost
$\nu$-stable derived equivalence between finite-dimensional algebras
$A$ and $B$  over a field $k$, then there is a stable equivalence
$\Phi$ of Morita type between $A$ and $B$ such that $\Phi(X)\simeq
\bar{F}(X)$ for all objects $X$ in the stable module category over
$A$. \label{thm1}
\end{Theo}

As a consequence of the proof of Theorem \ref{thm1}, we have the
following facts on the homological dimensions of algebras.

\begin{Koro} Let $A$ and $B$ be Artin algebras, and let $F$ be a derived
equivalence between $A$ and $B$. If $\add(_AQ)=\add(\nu_AQ)$, then

\smallskip
$(1)$ $\gldim (A)\leqslant \gldim(B)$,

\smallskip
$(2)$ $\findim(A)\leqslant \findim(B)$,

\smallskip
$(3)$  $\domdim(A)\geqslant \domdim(B)$.

\noindent where $\gldim(A)$, $\findim(A)$ and $\domdim(A)$ stand for
the global, finitistic and dominant dimensions of $A$, respectively.
\label{cor2}
\end{Koro}

Note that if $A$ and $B$ are finite-dimensional self-injective, we
re-obtain the well-known result \cite{RickDFun} of Rickard from
Theorem \ref{thm1}: Derived-equivalent self-injective algebras are
stably equivalent of Morita type. Moreover, Theorem \ref{thm1}
allows us to obtain a lot of (usually not self-injective) algebras
which are both derived-equivalent and stably equivalent of Morita
type. By the following corollary, we can even repeatedly construct
derived equivalences satisfying the almost $\nu$-stable condition.

\begin{Koro}
Let $k$ be a field, and let $F$ be an almost $\nu$-stable derived
equivalence between two finite-dimensional $k$-algebras $A$ and $B$.
Then:

$(1)$ For any finite-dimensional self-injective $k$-algebra $C$,
there is an almost $\nu$-stable derived equivalence between the two
tensor algebras $A\otimes_k C$ and $B\otimes_k C$.

$(2)$ Let $\bar{F}$ be the stable equivalence induced by $F$ in
Theorem \ref{thm1}. Then, for each $A$-module $X$, there is an
almost $\nu$-stable derived equivalence between the endomorphism
algebras $\End_A(A\oplus X)$ and $\End_B(B\oplus \bar{F}(X))$.

$(3)$ If $X$ is an $A$-module such that $F(X)$ is isomorphic to a
$B$-module $Y$, then there is an almost $\nu$-stable derived
equivalence between the one-point extensions $A[X]$ and $B[Y]$.
\end{Koro}

This paper is organized as follows. In Section \ref{preSection}, we
shall recall some basic definitions and facts required in proofs. In
Section \ref{theFunctor}, we first show that every derived
equivalence $F$ between two Artin algebras $A$ and $B$ gives rise to
a functor $\bar{F}$ between their stable module categories, and then
give a sufficient condition for the functor $\bar{F}$ to be an
equivalence between stable module categories over Artin algebras. In
Section \ref{dimensions}, we deduce some properties of the functor
$\bar{F}$ and then compare homological dimensions of $A$ with that
of $B$. In particular, we get Corollary \ref{cor2}. As a by-product,
we re-obtain the result that a derived equivalence preserves the
finiteness of finitistic dimension. In Section \ref{sectStM}, we
show that the condition given in Section \ref{theFunctor} is
sufficient for $F$ to induce a stable equivalence of Morita type
when we work with finite-dimensional algebras over a field. In
Section \ref{sectConst}, we give several methods to construct
inductively derived equivalences satisfying the almost $\nu$-stable
condition. Finally, in Section \ref{example}, we exhibit a couple of
examples to explain our points about the main result.

\section{Preliminaries}\label{preSection}
In this section, we shall recall basic definitions and facts
required in our proofs.

Let $\cal C$ be an additive category. For two morphisms
$f:X\rightarrow Y$ and $g:Y\rightarrow Z$ in $\cal C$, the
composition of $f$ with $g$ is written as $fg$, which is a morphism
from $X$ to $Z$. But for two functors $F:\mathcal{C}\rightarrow
\mathcal{D}$ and $G:\mathcal{D}\rightarrow\mathcal{E}$ of
categories, their composition is denoted by $GF$. For an object $X$
in $\mathcal{C}$, we denote by $\add(X)$ the full subcategory of
$\cal C$ consisting of all direct summands of finite direct sums of
copies of $X$.

Throughout this paper, unless specified otherwise, all algebras will
be Artin algebras over a fixed commutative Artin ring $R$. All
modules will be finitely generated unitary left modules. If $A$ is
an Artin algebra, the category of all modules over $A$ is denoted by
$\modcat{A}$; the full subcategory of $A$-mod consisting of
projective (respectively, injective) modules is denoted by
$\pmodcat{A}$ (respectively, $\imodcat{A}$). We  denote by $D$ the
usual duality on $\modcat{A}$. The duality $\Hom_A(-, A)$ from
$\pmodcat{A}$ to $\pmodcat{A\opp}$ is denoted by $^*$, that is, for
each projective $A$-module $P$, the projective $A\opp$-module
$\Hom_A(P, A)$ is denoted by $P^*$. We denote by $\nu_A$ the
Nakayama functor $D\Hom_A(-,A): A\mbox{-proj}\lra A\mbox{-inj}$.

The stable module category $\stmodcat{A}$ of an algebra $A$ is, by
definition, an $R$-category in which objects are the same as the
objects of $\modcat{A}$ and, for two objects $X, Y$ in
$\stmodcat{A}$, their morphism set, denoted by $\StHom_A(X,Y)$, is
the quotient of $\Hom_A(X,Y)$ modulo the homomorphisms that
factorize through projective modules. Two algebras are said to be
\emph{stably equivalent} if their stable module categories are
equivalent as $R$-categories.

For finite-dimensional algebras over a field $k$, there is a special
class of stable equivalences, namely stable equivalences of Morita
type. Recall that two finite-dimensional algebras $A$ and $B$ over a
field $k$ are called {\em stably equivalent of Morita type} if there
are two bimodules $_AM_B$ and $_BN_A$ satisfying the following
properties:

$(1)$ all of the one-sided modules $_AM, M_B, {}_BN$ and $N_A$  are
projective, and

$(2)$ there is an $A$-$A$-bimodule isomorphism
$_AM\otimes_BN_A\simeq A\oplus U$ for a projective $A$-$A$-bimodule
$U$; and there is a $B$-$B$-bimodule isomorphism
$_BN\otimes_AM_B\simeq B\oplus V$ for a projective $B$-$B$-bimodule
$V$.

\smallskip
If two finite-dimensional algebras $A$ and $B$ over a field are
stably equivalent of Morita type, then the functor $T_N:
\stmodcat{A}\lra\stmodcat{B}$ defined by $_BN\otimes_A-$ is an
equivalence. It is also called a {\em stable equivalence of Morita
type}. (Note that if we extend the definition of a stable
equivalence of Morita type from finite-dimensional algebras to
$R$-projective Artin $R$-algebras, then there is an open problem of
whether $T_N$ could induce a stable equivalence on stable module
categories, namely whether $_AU\otimes_A X$ is projective for every
$A$-module $X$.)

\medskip
Now, we recall some definitions relevant to derived equivalences.
Let $\cal C$ be an additive category. A complex $\cpx{X}$ over $\cal
C$ is a sequence of morphisms $d_X^{i}$ between objects $X^i$ in
$\cal C$: $ \cdots \rightarrow X^{i-1}\stackrel{d_X^{i-1}}{\lra}
X^i\stackrel{d_X^i}{\lra}
X^{i+1}\stackrel{d_X^{i+1}}{\lra}X^{i+2}\rightarrow\cdots $ such
that $d_X^id_X^{i+1}=0$ for all $i \in {\mathbb Z}$. We write
$\cpx{X}=(X^i, d_X^i)$. For each complex $\cpx{X}$, the {\em brutal
truncation} $\sigma_{<i}\cpx{X}$ is a subcomplex of $\cpx{X}$ such
that $(\sigma_{<i}\cpx{X})^k$ is $X^k$ for all $k<i$ and zero
otherwise. Similarly, we define $\sigma_{\geq i}\cpx{X}$. The
category of complexes over $\cal C$ is denoted by $\C{C}$. The
homotopy category of complexes over $\mathcal{C}$ is denoted by
$\K{\mathcal C}$. When $\cal C$ is an abelian category, the derived
category of complexes over $\cal C$ is denoted by $\D{\cal C}$. The
full subcategory of $\K{\cal C}$ and $\D{\cal C}$ consisting of
bounded complexes over $\mathcal{C}$ is denoted by $\Kb{\mathcal C}$
and $\Db{\mathcal C}$, respectively. Moreover, we denote by
$\Cf{\mathcal C}$ the category of complexes bounded above, and by
$\Kf{\mathcal C}$ the homotopy category of $\Cf{\cal C}$. Dually, we
have the category $\Cz{\cal C}$ of complexes bounded below and the
homotopy category $\Kz{\cal C}$ of $\Cz{\cal C}$. As usual, for a
given Artin algebra $A$, we simply write $\C{A}$ for
$\C{\modcat{A}}$, $\K{A}$ for $\K{\modcat{A}}$ and $\Kb{A}$ for
$\Kb{\modcat{A}}$. Similarly, we write $\D{A}$ and $\Db{A}$ for
$\D{\modcat{A}}$ and $\Db{\modcat{A}}$, respectively.

It is well-known that, for an Artin algebra $A$, $\K{A}$ and $\D{A}$
are triangulated categories. For basic results on triangulated
categories, we refer to Happel's book \cite{HappelTriangle}.
Throughout this paper, we denote by $X[n]$ rather than $T^nX$ the
object obtained from $X$ by shifting $n$ times. In particular, for a
complex $\cpx{X}$ in $\K{A}$ or $\D{A}$, the complex $\cpx{X}[1]$ is
obtained from $\cpx{X}$ by shifting $\cpx{X}$ to the left by one
degree.

Let $A$ be an Artin algebra. A homomorphism $f: X\lra Y$ of
$A$-modules is called a \emph{radical map} if, for any module $Z$
and homomorphisms $h: Z\lra X$ and $g: Y\lra Z$, the composition
$hfg$ is not an isomorphism. A complex over $\modcat{A}$ is called a
\emph{radical} complex if all of its differential maps are radical
maps. Every complex over $\modcat{A}$ is isomorphic in the homotopy
category $\K{A}$ to a radical complex. It is easy to see that if two
radical complex $\cpx{X}$ and $\cpx{Y}$  are isomorphic in $\K{A}$,
then $\cpx{X}$ and $\cpx{Y}$ are isomorphic in $\C{A}$.

Two algebras $A$ and $B$ are said to be \emph{derived-equivalent} if
their derived category $\Db{A}$ and $\Db{B}$ are equivalent as
triangulated categories. In \cite{RickMoritaTh}, Rickard proved that
two algebras are derived-equivalent if and only if there is a
complex $\cpx{T}$ in $\Kb{\pmodcat{A}}$ satisfying

(1) $\Hom_{\Db{A}}(\cpx{T},\cpx{T}[n])=0$ for all  $n\ne 0$, and

(2) $\add(\cpx{T})$ generates $\Kb{\pmodcat{A}}$ as a triangulated
category

{\parindent=0pt such that $B\simeq\End_{\Db{A}}(\cpx{T})$.} A
complex in $\Kb{\pmodcat{A}}$ satisfying the above two conditions is
called a {\em tilting complex} over $A$. By the condition (2), each
indecomposable projective $A$-module is a direct summand of $T^i$
for some integer $i$. It is known that, given a derived equivalence
$F$ between $A$ and $B$, there is a unique (up to isomorphism)
tilting complex $\cpx{T}$ over $A$ such that $F(\cpx{T})\simeq B$.
This complex $\cpx{T}$ is called a tilting complex
 \emph{associated} to $F$.

\medskip
Let $F$ be a derived equivalence between two Artin algebras $A$ and
$B$, and let $\cpx{Q}$ be a tilting complex associated to $F$.
Without loss of generality, we may assume that $\cpx{Q}$ is radical
and that $Q^i=0$ for $i<-n$ and $i>0$. Then we have the following
fact, for convenience of the reader, we include here a proof.

\begin{Lem}\label{tiltingforconormal}
Let $A$ and $B$ be two Artin algebras, and let $F$ and $\cpx{Q}$ be
as above. If $G:\Db{B}\lra\Db{A}$ is a quasi-inverse of $F$, then
there is a radical tilting complex $\cpx{\bar{Q}}$ associated to $G$
of the following form:
$$\xymatrix{0
\ar[r]& \bar{Q}^{0}\ar[r] & \bar{Q}^{1}\ar[r] &\cdots\ar[r] &
\bar{Q}^{n}\ar[r] & 0. }$$
\end{Lem}

{\it Proof.} Note that the tilting complex $\cpx{Q}$ associated to
$F$ is radical and of the form:
$$\xymatrix{0
\ar[r]& Q^{-n}\ar[r]&\cdots\ar[r]&Q^{-1}\ar[r] & Q^{0}\ar[r] & 0.}$$
Let $\cpx{\bar{Q}}$ be a radical complex in $\Kb{\pmodcat{B}}$ such
that $\cpx{\bar{Q}}$ isomorphic to $F(A)$ in $\Db{B}$. Then
$G(\cpx{\bar{Q}})\simeq GF(A)\simeq A$ in $\Db{A}$, which means that
$\cpx{\bar{Q}}$ is a tilting complex associated to $G$. Moreover, on
the one hand, we have
$$\Hom_{\Db{B}}(\cpx{\bar{Q}},B[m])\simeq \Hom_{\Db{A}}(A, \cpx{Q}[m])=0$$
for all $m>0$, and consequently $\cpx{\bar{Q}}$ has zero terms in
all negative degrees. On the other hand, we have
$$\Hom_{\Db{B}}(B,\cpx{\bar{Q}}[m])\simeq \Hom_{\Db{A}}(\cpx{Q}, A[m])=0$$
for all $m>n$ and $\Hom_{\Db{B}}(B,\cpx{\bar{Q}}[n])\neq 0$. Thus
the complex $\cpx{\bar{Q}}$ has zero terms in all degrees larger
than $n$, and its $n$-th term is non-zero. $\dickebox$

The following lemma will be  used frequently in our proofs. Again,
we include here a proof for convenience of the reader.

\begin{Lem}\label{kdiso}
Let $A$ be an arbitrary ring, and let $A\mbox{\rm -Mod}$ be the
category of all left (not necessarily finitely generated)
$A$-modules. Suppose $\cpx{X}$ is a complex over $A\mbox{\rm -Mod}$
bounded above and $\cpx{Y}$ is a complex over $A\mbox{\rm -Mod}$
bounded below. Let $m$ be an integer. If one of the following two
conditions holds:

$(1)$ $X^i$ is projective for all $i>m$ and $Y^j=0$ for all $j<m$;

$(2)$ $Y^j$ is injective for all $j<m$  and $X^i=0$ for all $i>m$,

\noindent then  $\theta_{\cpx{X},\cpx{Y}}: \Hom_{\K{A\mbox{\rm
-Mod}}}(\cpx{X},\cpx{Y})\rightarrow \Hom_{\D{A\mbox{\rm
-Mod}}}(\cpx{X},\cpx{Y})$ induced by the localization functor
$\theta: \K{A\mbox{\rm -Mod}}$ $\rightarrow\D{A\mbox{\rm -Mod}}$ is
an isomorphism.
\end{Lem}

{\it Proof.}  For simplicity, we write $\mathscr{K} = \K{A\mbox{\rm
-Mod}}$ and $\mathscr{D}=\D{A\mbox{\rm -Mod}}$. The category of all
left (not necessarily finitely generated) projective $A$-modules is
denoted by $A\mbox{-Proj}$. By applying the shift functor, we may
assume $m=0$. Suppose (1) is satisfied.

First, we consider the case $X^i=0$ for all $i<0$. Let
$$\xymatrix{ \cdots\ar[r] & P^{-1}\ar[r] &P^{0}\ar[r]^{\pi} & X^0\ar[r] & 0 }$$
be a projective resolution of $X^0$ with all $P^i$ being projective
$A$-modules. Then the following complex,
$$\xymatrix{\cdots\ar[r] & P^{-1}\ar[r] &P^{0}\ar[r]^{\pi d_X^0} & X^1 \ar[r] & X^2 \ar[r] & \cdots},$$
denoted by $\cpx{P_X}$, is in $\K{{A}\mbox{-Proj}}$ and bounded
above since $X^i$ is projective for all $i>0$ by our assumption, and
there is a quasi-isomorphism $\cpx{\pi}: \cpx{P_X}\rightarrow
\cpx{X}$:
$$\xymatrix{
\cdots\ar[r] & P^{-1}\ar[d]\ar[r] &P^{0}\ar[d]^{\pi}\ar[r]^{\pi
d_X^0} & X^1\ar@{=}[d] \ar[r]
&\cdots\\
&0\ar[r]& X^0\ar[r]^{d_X^0} & X^1\ar[r] &\cdots. }$$ We claim that
${\Hom_{\mathscr{K}}(\cpx{\pi}, \cpx{Y})}:
\Hom_{\mathscr{K}}(\cpx{X},\cpx{Y})\lra
\Hom_{\mathscr{K}}(\cpx{P_X},\cpx{Y})$ induced by $\cpx{\pi}$ is an
isomorphism. Actually, if $\cpx{f}\in\Hom_{\mathscr{K}}(\cpx{P_X},
\cpx{Y})$, then $f^0$ factorizes through the map $\pi:
P^0\rightarrow X^0$. Suppose $f^0=\pi g^0$ for some $g^0: X^0\ra
Y^0$. Let $g^i := f^i$ for all $i>0$. Then $\cpx{g}=(g^i)$ is a
chain map from $\cpx{X}$ to $\cpx{Y}$ and
$\cpx{f}=\cpx{\pi}\cpx{g}$. Consequently, the map
$\Hom_{\mathscr{K}}(\cpx{\pi}, \cpx{Y})$ is surjective. Further, we
show that the map $\Hom_{\mathscr{K}}(\cpx{\pi}, \cpx{Y})$ is
injective. In fact, if $\cpx{\pi}\cpx{\alpha}=0$ for a morphism
$\cpx{\alpha}$ in $\Hom_{\mathscr{K}}(\cpx{X},\cpx{Y})$, then there
are maps $h^i: X^i\rightarrow Y^{i-1}$ for all integer $i\ge 1$ such
that $\pi\alpha^0=\pi d_X^0h^1$ and
$\alpha^i=d_X^ih^{i+1}+h^id_Y^{i-1}$. Thus $\alpha^0=d_X^0h^1$ since
$\pi$ is an epimorphism. Hence $\cpx{\alpha}=0$, which implies that
$\Hom_{\mathscr{K}}(\cpx{\pi}, \cpx{Y})$ is injective. It follows
that $\Hom_{\mathscr{K}}(\cpx{\pi}, \cpx{Y})$ is an isomorphism.

Note that $\cpx{\pi}$ induces a commutative diagram
$$\begin{CD}
\Hom_{\mathscr{K}}(\cpx{X},\cpx{Y})@>{\theta_{\cpx{X},\cpx{Y}}}>>
\Hom_{\mathscr{D}}(\cpx{X},\cpx{Y})\\
@V{(\cpx{\pi}, -)}VV @V{(\cpx{\pi}, -)}VV\\
\Hom_{\mathscr{K}}(\cpx{P_X},\cpx{Y})@>{\theta_{\cpx{P_X},\cpx{Y}}}>>
\Hom_{\mathscr{D}}(\cpx{P_X},\cpx{Y}).\\
\end{CD}$$
Since $\cpx{\pi}$ is a quasi-isomorphism, the right vertical map is
an isomorphism. We have shown that the left vertical map is an
isomorphism. Moreover, since $\cpx{P_X}$ is a complex in
$\K{A\mbox{-Proj}}$ and bounded above, the map
$\theta_{\cpx{P_X},\cpx{Y}}$ is an isomorphism. It follows that
$\theta_{\cpx{X},\cpx{Y}}:
\Hom_{\mathscr{K}}(\cpx{X},\cpx{Y})\longrightarrow
\Hom_{\mathscr{D}}(\cpx{X},\cpx{Y})$ is an isomorphism.

Now, let $\cpx{X}$ be an arbitrary complex satisfying the condition
(1). Then there is a distinguished triangle
$$\sigma_{<0}\cpx{X}[-1]\longrightarrow\sigma_{\geq 0}\cpx{X}
\stackrel{p}{\longrightarrow}\cpx{X}\longrightarrow\sigma_{<0}\cpx{X}$$
in $\mathscr{K}$. This triangle can be viewed as a distinguished
triangle in $\mathscr{D}$. Applying the functors
$\Hom_{\mathscr{K}}(-, \cpx{Y})$ and $\Hom_{\mathscr{D}}(-,
\cpx{Y})$ to the triangle, we have an exact commutative diagram
$$\begin{CD}
\Hom_{\mathscr{K}}(\sigma_{<0}\cpx{X}, \cpx{Y})@>>>
\Hom_{\mathscr{K}}(\cpx{X}, \cpx{Y})@>{{}_{\mathscr{K}}(p,-)}>>
\Hom_{\mathscr{K}}(\sigma_{\geq 0}\cpx{X}, \cpx{Y})
@>>>\Hom_{\mathscr{K}}(\sigma_{<0}\cpx{X}[-1], \cpx{Y})\\
@VVV @V{\theta_{\cpx{X},\cpx{Y}}}VV @V{\theta_{\sigma_{\geq
0}\cpx{X},\cpx{Y}}}VV
@V{\theta_{\sigma_{< 0}\cpx{X}[-1],\cpx{Y}}}VV\\
\Hom_{\mathscr{D}}(\sigma_{<0}\cpx{X}, \cpx{Y})@>>>
\Hom_{\mathscr{D}}(\cpx{X}, \cpx{Y})@>{{}_{\mathscr{D}}(p,-)}>>
\Hom_{\mathscr{D}}(\sigma_{\geq 0}\cpx{X}, \cpx{Y})@>>>
\Hom_{\mathscr{D}}(\sigma_{<0}\cpx{X}[-1], \cpx{Y}).\\
\end{CD}$$
By our assumption, we have $Y^i=0$ for all $i<0$, and therefore
$\Hom_{\mathscr{K}}(\sigma_{<0}\cpx{X}, \cpx{Y})=0$. Note that
$\sigma_{<0}\cpx{X}$ is isomorphic in $\mathscr{D}$ to a complex
$\cpx{P_1}$ in $\K{{A}\mbox{-Proj}}$ such that $P_1^i=0$ for all
$i\geq 0$. It follows that $\Hom_{\mathscr{D}}(\sigma_{<0}\cpx{X},
\cpx{Y})\simeq \Hom_{\mathscr{D}}(\cpx{P_1}, \cpx{Y})\simeq
\Hom_{\mathscr{K}}(\cpx{P_1}, \cpx{Y})=0$. Thus both maps
${{}_{\mathscr{K}}(p,-)}$ and ${{}_{\mathscr{D}}(p,-)}$ are
injective. Note that $\sigma_{\geq 0}\cpx{X}$ is a complex
satisfying the condition (1) and the terms of $\sigma_{\geq
0}\cpx{X}$ in all negative degrees are zero. By the first part of
the proof, we see that the map $\theta_{\sigma_{\geq
0}\cpx{X},\cpx{Y}}$ is an isomorphism. It follows that
$\theta_{\cpx{X},\cpx{Y}}$ is injective. In particular, the map
$\theta_{\sigma_{< 0}\cpx{X}[-1],\cpx{Y}}$ is injective. By the Five
Lemma (see \cite[Exercise 1.3.3, p.13]{Weibel}, the map
$\theta_{\cpx{X},\cpx{Y}}$ is surjective. Thus
$\theta_{\cpx{X},\cpx{Y}}$ is an isomorphism.

The proof for the situation (2) proceeds dually. $\dickebox$

\smallskip
In the following, we point out a relationship between the Nakayama
functor and a derived equivalence. Let $F:\Db{A}\lra\Db{B}$ be a
derived equivalence between Artin algebras $A$ and $B$ over a
commutative Artin ring $R$. By \cite[Theorem 6.4]{RickMoritaTh}, $F$
induces an equivalence from $\Kb{\pmodcat{A}}$ to
$\Kb{\pmodcat{B}}$. Note that the Nakayama functor $\nu_A:
\pmodcat{A}\lra \imodcat{A}$ induces a functor from
$\Kb{\pmodcat{A}}$ to $\Kb{\imodcat{A}}$, which is again denoted by
$\nu_A$. When $A$ and $B$ are finite-dimensional algebras over a
field $k$, it is known that $F(\nu_A\cpx{P})\simeq\nu_BF(\cpx{P})$
in $\Db{B}$ for each $\cpx{P}\in\Kb{\pmodcat{A}}$ \cite{RickDFun}.
We shall extend this to Artin algebras by using the notation of
Auslander-Reiten triangle. Recall that a distinguished triangle
$X\stackrel{f}{\lra}M\stackrel{g}{\lra}Y\stackrel{w}{\lra}X[1]$ in
$\Db{A}$ is called an {\em Auslander-Reiten triangle} if

(AR1) $X$ and $Y$ are indecomposable,

(AR2) $w\ne 0$, and

(AR3) if $t:U\longrightarrow Y$ is not a split epimorphism, then
$tw=0$.

For a given object $Y$ in $\Db{A}$, if there is an Auslander-Reiten
triangle
$X\stackrel{f}{\lra}M\stackrel{g}{\lra}Y\stackrel{w}{\lra}X[1]$,
then it is unique up to isomorphism \cite[Proposition 4.3,
p.33]{HappelTriangle}.

The first statement of the following lemma is essentially due to
Happel (see \cite[Theorem 4.6, p.37]{HappelTriangle}), here we shall
modify the original proof to deal with Artin algebras.

\begin{Lem}
Let $A$ be an Artin algebra over a commutative Artin ring $R$. Then,
for each indecomposable object $\cpx{P}$ in $\Kb{\pmodcat{A}}$,
there is an Auslander-Reiten triangle
$$(\nu_{A}\cpx{P})[-1]\longrightarrow \cpx{L}\longrightarrow
\cpx{P}\longrightarrow\nu_{A}\cpx{P}$$ in $\Db{A}$. Furthermore, we
have $ F(\nu_A\cpx{P}))\simeq \nu_BF(\cpx{P})$ in $\Db{B}$.
\label{nakayama}
\end{Lem}

{\it Proof.} Note that there is an invertible natural transformation
$\alpha_P: D\Hom_{A}(P, -)\longrightarrow \Hom_{A}(-,\nu_{A}P)$ for
each projective $A$-module $P$. This induces an invertible natural
transformation of two functors from $\Db{A}$ to $R\modcat$.
$$\alpha_{\cpx{P}}: D\Hom_{\Db{A}}(\cpx{P}, -)\longrightarrow
\Hom_{\Db{A}}(-,\nu_{A}\cpx{P}).$$ Recall that $D=\Hom_R(-, J)$ with
$J$ an injective envelope of the $R$-module $R/\mbox{rad}(R)$. Let
$\psi$ be a non-zero $R$-module homomorphism from
$\End_{\Db{A}}(\cpx{P})/\rad(\End_{\Db{A}}(\cpx{P}))$ to $J$. We
define $\phi: \End_{\Db{A}}(\cpx{P})\lra J$ to be the composition of
the canonical surjective homomorphism $$\End_{\Db{A}}(\cpx{P})\lra
\End_{\Db{A}}(\cpx{P})/\rad(\End_{\Db{A}}(\cpx{P}))$$ with $\psi$ .
Then $\phi$ is non-zero and vanishes on
$\rad(\End_{\Db{A}}(\cpx{P}))$. Thus we have a distinguished
triangle
$$\begin{CD}
(\nu_{A}\cpx{P})[-1]@>>> \cpx{L}@>>>
\cpx{P}@>{\alpha_{\cpx{P}}(\phi)}>>\nu_{A}\cpx{P},
\end{CD}$$ where $\cpx{L}[1]$ is the mapping cone of the map $\alpha_{\cpx{P}}(\phi)$.
Clearly, this triangle satisfies the conditions (AR1) and (AR2). Let
$f: \cpx{X}\longrightarrow\cpx{P}$ be a morphism which is not a
split epimorphism. Then we have a commutative diagram
$$\begin{CD}
D\Hom_{\Db{A}}(\cpx{P}, \cpx{P}) @>{\alpha_{\cpx{P}}}>>
\Hom_{\Db{A}}(\cpx{P},
\nu_{A}\cpx{P})\\
@V(f_*, -)VV @VV(f,-)V\\
D\Hom_{\Db{A}}(\cpx{P}, \cpx{X}) @>{\alpha_{\cpx{P}}}>>
\Hom_{\Db{A}}(\cpx{X},
\nu_{A}\cpx{P}),\\
\end{CD}$$
where $f_*: \Hom_{\Db{A}}(\cpx{P}, \cpx{P})\lra
\Hom_{\Db{A}}(\cpx{P}, \cpx{X})$ is induced by $f$. Since $f$ is not
a split epimorphism and $\cpx{P}$ is indecomposable, we  find that,
for each morphism $g: \cpx{P}\longrightarrow\cpx{X}$, the
composition $gf$ is in $\rad(\End_{\Db{A}}(\cpx{P}))$. By the
definition of $\phi$, we have $(f_*\phi)(g)=\phi(gf)=0$ for all
$g\in\Hom_{\Db{A}}(\cpx{P}, \cpx{X})$, that is, $f_*\phi=0$. It
follows that $f\alpha_{\cpx{P}}(\phi)=\alpha_{\cpx{P}}(f_*\phi)=0$.
This proves (AR3).

Now, for the indecomposable object $\cpx{P}$ in $\Kb{\pmodcat{A}}$,
there is an Auslander-Reiten triangle
$$(*)\quad (\nu_{A}\cpx{P})[-1]\longrightarrow \cpx{L}\longrightarrow
\cpx{P}\longrightarrow\nu_{A}\cpx{P} $$ in $\Db{A}$. Since
$F(\cpx{P})$ is in $\Kb{\pmodcat{B}}$ and indecomposable, there is
also an Auslander-Reiten triangle
$$(\nu_{B}F(\cpx{P}))[-1]\longrightarrow \cpx{L'}\longrightarrow
F(\cpx{P})\longrightarrow\nu_{B}F(\cpx{P}) $$ in $\Db{B}$. Further,
if we apply the functor $F$ to ($*$), we get another
Auslander-Reiten triangle
$$F(\nu_{A}\cpx{P})[-1]\longrightarrow F(\cpx{L})\longrightarrow
F(\cpx{P})\longrightarrow F(\nu_{A}\cpx{P}) $$ in $\Db{B}$.  By the
uniqueness of the Auslander-Reiten triangle associated to the given
complex $F(\cpx{P})$, we see that $\nu_BF(\cpx{P})$ is isomorphic to
$F(\nu_{A}\cpx{P})$ in $\Db{B}$. $\dickebox$

\medskip
Finally, let us remark that, given a functor
$F:\mathcal{C}\rightarrow\mathcal{D}$, if we fix an object $F_X$ in
$\cal D$ for each object $X$ in $\mathcal{C}$ such that $F_X\simeq
F(X)$, then there is a functor
$F':\mathcal{C}\rightarrow\mathcal{D}$ such that $F'\simeq F$ and
$F'(X)=F_X$ for every $X$ in $\mathcal{C}$. Actually, let $s_X$
denote the isomorphism from $F_X$ to $F(X)$, and we define
$F'(f):=s_XF(f)s_Y^{-1}$  for each $f: X\rightarrow Y$. Then this
$F'$ is a desired functor.

\section{Stable equivalences induced by derived equivalences}\label{theFunctor}

In this section, we shall first construct a functor $\bar{F}:
\stmodcat{A}\lra\stmodcat{B}$ between the stable module categories
of two Artin algebras $A$ and $B$ from a given derived equivalence
$F:\Db{A}\lra\Db{B}$, and then give a sufficient condition to ensure
that the functor $\bar{F}$ is an equivalence. In Section
\ref{sectStM}, we shall see a stronger conclusion when we work with
finite-dimensional algebras instead of general Artin algebras.

Let us first recall some notions and notations. Let $A$ be an Artin
algebra. The homotopy category $\Kb{\pmodcat{A}}$ can be considered
as a triangulated full subcategory of $\Db{A}$. Let
$\Db{A}/\Kb{\pmodcat{A}}$ be the Verdier quotient of $\Db{A}$ by the
subcategory $\Kb{\pmodcat{A}}$ (for the definition, we refer the
reader to the excellent book \cite{neumann}). Then there is a
canonical functor $\Sigma': \modcat{A}\longrightarrow
\Db{A}/\Kb{\pmodcat{A}}$ obtained by composing the natural embedding
of $\modcat{A}$ into $\Db{A}$ with the quotient functor from
$\Db{A}$ to $\Db{A}/\Kb{\pmodcat{A}}$. Clearly, $\Sigma'(P)$ is
isomorphic to zero for each projective $A$-module $P$, so $\Sigma'$
factorizes through the natural functor $\modcat{A}\lra\stmodcat{A}$.
This gives rise to an additive functor $\Sigma:
\stmodcat{A}\lra\Db{A}/\Kb{\pmodcat{A}}$.

Rickard \cite{RickDstable} proved that $\Sigma$ is an equivalence
provided that the algebra $A$ is self-injective. But for an
arbitrary algebra, this is no longer true in general; for instance,
if $A$ is a non-semisimple Artin algebra of finite global dimension,
then the quotient category $\Db{A}/\Kb{\pmodcat{A}}$ is zero, and
therefore the functor $\Sigma$ is a zero functor which cannot be an
equivalence.

Let $A$ and $B$ be two Artin algebras. Suppose $F:
\Db{A}\longrightarrow\Db{B}$ is a derived equivalence between $A$
and $B$. Then $F$ induces an equivalence between the quotient
categories $\Db{A}/\Kb{\pmodcat{A}}$ and $\Db{B}/\Kb{\pmodcat{B}}$.
For simplicity, we denote this induced equivalence also by $F$.
Thus, if $A$ and $B$ are self-injective, then $\stmodcat{A}$ and
$\stmodcat{B}$ are equivalent. However, this is not true in general
for arbitrary finite-dimensional algebras, namely we cannot get an
equivalence of stable module categories from a given derived
equivalence in general. Nevertheless, we may ask if there is any
$``$good" functor $\bar{F}: \stmodcat{A}\lra \stmodcat{B}$ induced
by $F$, which could be a possible candidate for a stable equivalence
under certain additional conditions, and would cover the most
interesting known situations.

In the following, we shall construct an additive functor $\bar{F}:
\stmodcat{A}\lra\stmodcat{B}$ from $F$ such that the diagram
$$\begin{CD}
\stmodcat{A}@>{\Sigma}>>\Db{A}/\Kb{\pmodcat{A}}\\
@V{\bar{F}}VV @VVFV\\
 \stmodcat{B}@>\Sigma>> \Db{B}/\Kb{\pmodcat{B}}
\end{CD}$$
of additive functors is commutative up to natural isomorphisms.
Furthermore, we shall construct a possible candidate for the inverse
of $\bar{F}$ under an additional condition.

From now on, $A$ and $B$ are  Artin $R$-algebras, $F$ is a derived
equivalence between $A$ and $B$ with the quasi-inverse $G$. Let
$\cpx{Q}$ be a tilting complex over $A$ associated to $F$ of the
following form:
      $$\xymatrix{ \cpx{Q}: \qquad
      0 \ar[r]& Q^{-n}\ar[r]&\cdots\ar[r]&Q^{-1}\ar[r] & Q^0\ar[r] & 0
      }$$
such that all differentials are radical maps. By Lemma
\ref{tiltingforconormal}, there is a tilting complex $\cpx{\bar{Q}}$
associated to $G$ of the form
$$\xymatrix{ \cpx{\bar{Q}}: \qquad
0\ar[r]& \bar{Q}^0\ar[r] & \bar{Q}^{1}\ar[r] &\cdots\ar[r] &
\bar{Q}^{n}\ar[r] & 0 }$$ with all differentials being radical maps.
We define $Q=\bigoplus_{i=1}^{n}Q^{-i}$ and
$\bar{Q}=\bigoplus_{i=1}^{n}\bar{Q}^i$.

\begin{Lem}\label{formofFX}
Let $X$ be an $A$-module. Then $F(X)$ is isomorphic in $\Db{B}$ to a
radical complex $\cpx{\bar{Q}_X}$ of the following form
$$\xymatrix{0
\ar[r]& \bar{Q}_X^0\ar[r] & \bar{Q}_X^{1}\ar[r] &\cdots\ar[r] &
\bar{Q}_X^{n}\ar[r] & 0 },$$ with $\bar{Q}_X^i\in\add(_B\bar{Q})$
for all $i=1,2,\cdots,n$. Moreover, the complex $\cpx{\bar{Q}_X}$ of
this form is unique up to isomorphisms in $\Cb{B}$. In particular,
if $X$ is projective, then $\cpx{\bar{Q}}_X$ is isomorphic in
$\Cb{B}$ to a complex in $\add(\cpx{\bar{Q}})$. \label{image}
\end{Lem}

{\it Proof.} Let $H^i$ be the $i$-th homology functor on complexes.
First of all, we have $H^i(F(X)) \simeq \Hom_{\Db{B}}(B,F(X)[i])$
$\simeq \Hom_{\Db{A}}(\cpx{Q},X[i]) = 0$ for all $i>n$ and all
$i<0$, which means that $F(X)$ has no homology in negative degrees
and degrees larger than $n$. Clearly, we may assume that $X$ is
indecomposable.

If $X$ is projective, then $X$ is isomorphic to a direct summand of
$A$. Consequently, $F(X)$ is isomorphic in $\Db{B}$ to a direct
summand $\cpx{L}$ of the complex $\cpx{\bar{Q}}$. Since all terms of
$\cpx{\bar{Q}}$ in positive degrees are in $\add(_B\bar{Q})$, all
terms of $\cpx{L}$ in positive degrees are in $\add(_B\bar{Q})$.
This shows that for every projective $A$-module $P$, the complex
$F(P)$ is isomorphic in $\Db{B}$ to a complex with all of its terms
of positive degrees in $\add(_B\bar{Q})$. Now we show that if
$\cpx{P}$ is a complex in $\Kb{\pmodcat{A}}$ with $P^i=0$ for all
$i>0$, then $F(\cpx{P})$ is isomorphic in $\Db{B}$ to a complex in
which all of its terms in positive degrees belong to
$\add(_B\bar{Q})$. In fact, if $\cpx{P}$ has only one non-zero term,
then we may write $\cpx{P}=P[t]$ for a projective $A$-module $P$ and
a non-negative integer $t$. In this case, $F(\cpx{P})$ is isomorphic
to a direct summand of $\cpx{\bar{Q}}[t]$ in which  all terms in
positive degrees are in $\add(_B\bar{Q})$, as desired. Now, we
assume that $\cpx{P}$ has at least two non-zero terms. Then there is
an integer $s<0$ such that the brutal truncations
$\sigma_{<s}\cpx{P}$ and $\sigma_{\geq s}\cpx{P}$ have less non-zero
terms than $\cpx{P}$ does. By induction, the complexes
$F(\sigma_{<s}\cpx{P})$ and $F(\sigma_{\geq s}\cpx{P})$ are
respectively isomorphic to complexes $\cpx{Y}$ and $\cpx{Z}$ in
$\Kb{\pmodcat{A}}$, such that their terms in all positive degrees
are in $\add(_B\bar{Q})$. Since $\cpx{P}$ is the mapping cone of the
map $d_P^{s-1}$ from $\sigma_{<s}\cpx{P}$ to $\sigma_{\geq
s}\cpx{P}$, the complex $F(\cpx{P})$ is isomorphic to the mapping
cone of a chain map from $\cpx{Y}$ to $\cpx{Z}$, and consequently
all of its terms in positive degrees lie in $\add(_B\bar{Q})$.

Now, suppose that $X$ is an arbitrary indecomposable $A$-module and
$\cpx{P}=(P^i, d^i)$ is a minimal projective resolution of $X$. We
denote by $\Omega^{n}(X)$ the $n$-th syzygy of $X$, and by
$\cpx{P_1}$  the complex
$$\xymatrix{0
\ar[r]& P^{-n+1}\ar[r] & P^{-n+2}\ar[r] &\cdots\ar[r] & P^0 \ar[r] &
0. }$$ Then we have a distinguished triangle in $\Db{A}$
$$\xymatrix{
\Omega^{n}(X)[n-1]\ar[r] & \cpx{P_1}\ar[r] & X\ar[r] &
\Omega^{n}(X)[n].}$$  From this triangle one gets the following
distinguished triangle in $\Db{B}$:
$$\xymatrix{F(\Omega^{n}(X))[n-1]\ar[r] & F(\cpx{P_1})\ar[r] & F(X)\ar[r] &
F(\Omega^{n}(X))[n]. }$$ The complex $\cpx{P_1}$ is in
$\Kb{\pmodcat{A}}$ and all the terms of $\cpx{P_1}$ in positive
degrees are zero. Hence $F(\cpx{P_1})$ is isomorphic to a complex
$\cpx{Q_1}$ in $\Kb{\pmodcat{B}}$ with $Q_1^i$ in $\add(_B\bar{Q})$
for all $i>0$. Since $\Omega^{n}(X)$ is an $A$-module, the complex
$F(\Omega^{n}(X))$ has no homology in all degrees larger than $n$.
Thus the complex $F(\Omega^{n}(X))$ is isomorphic in $\D{B}$ to a
complex $\cpx{P_2}\in\Kf{\pmodcat{B}}$ with zero terms in all
degrees larger than $n$. It follows that $\cpx{P_2}[n-1]$ has zero
terms in all degrees larger than $1$. Hence $F(X)$ is isomorphic to
the mapping cone $\mbox{con}(\mu)$ of a map $\mu$ from
$\cpx{P_2}[n-1]$ to $\cpx{Q_1}$, and all the terms of
$\mbox{con}(\mu)$ in positive degrees are in $\add(_B\bar{Q})$. Note
that $F(X)$ has zero homology in all negative degrees and degrees
larger than $n$. Thus $\mbox{con}(\mu)$ has the same property. Hence
$\mbox{con}(\mu)$ is isomorphic in $\D{B}$ to a radical complex
$$\xymatrix{0
\ar[r]& \bar{Q}_X^0\ar[r] & \bar{Q}_X^{1}\ar[r] &\cdots\ar[r] &
\bar{Q}_X^{n}\ar[r] & 0 }$$ with $\bar{Q}_X^i\in\add(_B\bar{Q})$ for
all $i = 1, 2, \cdots, n$.

Suppose $\cpx{U}$ and $\cpx{V}$ are two radical complexes of the
form in Lemma \ref{image} such that both $\cpx{U}$ and $\cpx{V}$ are
isomorphic to $F(X)$ in $\Db{B}$. Then $\cpx{U}$ and $\cpx{V}$ are
isomorphic in $\Kb{B}$ by Lemma \ref{kdiso}. Since $\cpx{U}$ and
$\cpx{V}$ are radical complexes, we know that $\cpx{U}$ and
$\cpx{V}$ are isomorphic as complexes.

If $X$ is projective, then $X\in \add(A)$ and $F(X)\in \add(F(A))$.
Since $F(A)\simeq FG(\cpx{\bar{Q}})\simeq \cpx{\bar{Q}}$ in
$\Db{B}$, we see that $F(X)$ is isomorphic in $\Kb{\pmodcat{B}}$ to
a complex $\cpx{Y}\in \add(\cpx{\bar{Q}})$.  Thus $\cpx{\bar{Q}}_X$
is isomorphic in $\Db{B}$ to the $\cpx{Y}$. Since $\cpx{Y}$ is a
complex with the properties in Lemma \ref{image}, we have
$\cpx{Y}\simeq \cpx{\bar{Q}}_X$ as complexes by the uniqueness of
$\cpx{\bar{Q}}_X$. This shows that $\cpx{\bar{Q}}_X\in
\add(\cpx{\bar{Q}})$. Thus Lemma \ref{image} is proved. $\dickebox$

\medskip
Dually, we have the following lemma.
\begin{Lem}\label{formofGY}
Let $Y$ be a $B$-module. Then $G(Y)$ is isomorphic in $\Db{A}$ to a
radical complex $\cpx{Q_Y}$ of the form
$$\xymatrix{
0 \ar[r]&  Q_Y^{-n}\ar[r]&\cdots\ar[r]& Q_Y^{-1}\ar[r] & Q_Y^0\ar[r]
& 0 }$$ with $Q_Y^{-i}\in\add (\nu_AQ)$ for all $i=1, 2, \cdots, n$.
Moreover, the complex $\cpx{Q_Y}$ of this form is unique up to
isomorphisms in $\Cb{A}$.
\end{Lem}

{\it Remark. }One can easily see that if $X\simeq Y\oplus Z$ in
$\modcat{A}$, then the complex $\cpx{\bar{Q}_{X}}$ defined in Lemma
\ref{formofFX} is isomorphic in $\Cb{B}$ to the direct sum of
$\cpx{\bar{Q}_{Y}}$ and $\cpx{\bar{Q}_{Z}}$. Similarly, if $U\simeq
V\oplus W$ in $\modcat{B}$, then the complex $\cpx{Q_U}$ defined in
Lemma \ref{formofGY} is isomorphic in $\Cb{A}$ to the direct sum of
$\cpx{Q_V}$ and $\cpx{Q_W}$.

\smallskip
The next lemma is useful in our proofs.

\begin{Lem}\label{factthroughP}
Let $A$ be an Artin algebra, and let $f: X\lra Y$ be a homomorphism
between two $A$-modules $X$ and $Y$. Suppose $\cpx{P}$ is a complex
in $\Kb{A}$ with $P^i$ projective for all $i\geq 0$ and injective
for all $i<0$. If $f$ factorizes in $\Db{A}$ through $\cpx{P}$, then
$f$ factorizes through a projective $A$-module.
\end{Lem}

{\it Proof.}  There is a distinguished triangle
$$  \sigma_{<1}\cpx{P}[-1]\lra \sigma_{\geq
1}\cpx{P}\stackrel{b}{\lra}
\cpx{P}\stackrel{a}{\lra}\sigma_{<1}\cpx{P}$$ in $\Db{A}$. Note that
$\Hom_{\Db{A}}(\sigma_{\geq 1}\cpx{P}, Y)\simeq
\Hom_{\Kb{A}}(\sigma_{\geq 1}\cpx{P}, Y)=0$. Thus, if $f=gh$ for a
morphism $g: X\lra \cpx{P}$ and a morphism $h: \cpx{P}\lra Y$, then
$bh=0$, and consequently $h$ factorizes through $a$, say $h=ah'$ for
$h': \sigma_{<1}\cpx{P}\lra Y$. This means that $f=gah'$ and
factorizes through $\sigma_{<1}\cpx{P}$. Thus we may assume $P^i=0$
for all $i>0$ and consider the following distinguished triangle
$$\xymatrix{
  \sigma_{<0}\cpx{P}[-1]\ar[r] &  P^0\ar[r]^{u} & \cpx{P}\ar[r]^{v} & \sigma_{<0}\cpx{P}
 }$$
in $\Db{A}$. Now, we suppose $f=gh$ for $g: X\lra \cpx{P}$ and
$h:\cpx{P}\lra Y$. Note that the complex $\sigma_{<0}\cpx{P}$ is in
$\Kb{\imodcat{A}}$ by our assumption. Hence $\Hom_{\Db{A}}(X,
\sigma_{<0}\cpx{P})\simeq \Hom_{\Kb{A}}(X, \sigma_{<0}\cpx{P})=0$,
and consequently $gv=0$. This implies that $g$ factorizes through
$P^0$, that is, $g=g'u$ for a morphism $g': X\longrightarrow P^0$.
Since $\modcat{A}$ is fully embedded in $\Db{A}$, the morphisms $g'$
and $uh$ are $A$-module homomorphisms, and therefore $f=gh=g'(uh)$,
which factorizes through the projective $A$-module $P^0$.
$\dickebox$

\medskip
Now we define the functor $\bar{F}$. Pick an $A$-module $X$, by
Lemma \ref{formofFX}, we know that $F(X)$ is isomorphic in $\Db{B}$
to a radical complex $\cpx{\bar{Q}_X}$ of the form
$$\xymatrix{0
\ar[r]& \bar{Q}_X^0\ar[r] & \bar{Q}_X^{1}\ar[r] &\cdots\ar[r] &
\bar{Q}_X^{n}\ar[r] & 0 }$$ with $\bar{Q}_X^i\in\add(_B\bar{Q})$ for
all $i=1,2,\cdots,n$. From now on, for each $A$-module $X$, we
choose (once and for all) such a complex $\cpx{\bar{Q}_X}$. For each
homomorphism $f: X\longrightarrow Y$, we denote by $\underline{f}$
the image of $f$ under the canonical surjective map from $\Hom_A(X,
Y)$ to $\StHom_A(X,Y)$.

\begin{Prop}\label{stablefunctor}
Let $F:\Db{A}\lra\Db{B}$ be a derived equivalence between Artin
algebras $A$ and $B$. Then there is an additive functor $\bar{F}:
\stmodcat{A}\lra\stmodcat{B}$ sending $X$ to $\bar{Q}_X^0$ such that
the following diagram of the additive functors
$$\begin{CD}
\stmodcat{A}@>{\Sigma}>>\Db{A}/\Kb{\pmodcat{A}}\\
@V{\bar{F}}VV @VVFV\\
 \stmodcat{B}@>{\Sigma}>> \Db{B}/\Kb{\pmodcat{B}}
\end{CD}$$
is commutative up to natural isomorphisms.
\end{Prop}

{\it Proof.} By the remark at the end of Section 2 and Lemma
\ref{formofFX}, we may assume that $F(X)$ is just $\cpx{\bar{Q}_X}$
for each $A$-module $X$, where $\cpx{\bar{Q}_X}$ is the complex that
we have fixed above. Let $\bar{Q}_X^+$ denote the complex
$\sigma_{\geq 1}\cpx{\bar{Q}_X}$. Then we have a distinguished
triangle in $\Db{B}$:
$$\xymatrix{
\bar{Q}_X^+\ar[r]^{i_X} & F(X)\ar[r]^{\pi_X} &
\bar{Q}_X^0\ar[r]^(.40){\alpha_X} & \bar{Q}_X^+[1]. }$$ For each
homomorphism $f: X\lra Y$ of $A$-modules $X$ and $Y$, there is a
commutative diagram
$$\begin{CD}  \bar{Q}_X^+@>{i_X}>>F(X)@>{\pi_X}>>\bar{Q}_X^0 @>{\alpha_X}>>\bar{Q}_X^+[1]\\
@VV{a_f}V @VV{F(f)}V @VV{b_f}V @VV{a_f[1]}V \\
\bar{Q}_Y^+ @>{i_Y}>>F(Y)@>{\pi_Y}>>\bar{Q}_Y^0
@>{\alpha_Y}>>\bar{Q}_Y^+[1].
\end{CD}$$
The map $a_f$ exists because $i_XF(f)\pi_Y$ belongs to
$\Hom_{\Db{B}}(\bar{Q}_X^+,Q^0_Y)\simeq
\Hom_{\Kb{B}}(\bar{Q}_X^+,Q^0_Y)=0$. Since $B$-mod is fully embedded
in $\D{B}$, the morphism $b_f$ is a homomorphism of modules. If we
have another $A$-module homomorphism $b_f'$ such that
$\pi_Xb_f'=F(f)\pi_Y$, then $\pi_X(b_f-b_f')=F(f)\pi_Y-F(f)\pi_Y=0$
and $b_f-b_f'$ factorizes through $\bar{Q}_X^+[1]$. By Lemma
\ref{factthroughP}, the $B$-module homomorphism $b_f-b_f'$
factorizes through a projective $B$-module. Thus, for each
$A$-module homomorphism $f$ in $\Hom_A(X, Y)$, the morphism
$\underline{b_f}$ in  $\StHom_B(\bar{Q}_X^0, \bar{Q}_Y^0)$ is
uniquely determined by $f$.

\smallskip
Suppose $f: X\lra Y$ and $g: Y\lra Z$ are two homomorphisms of
$A$-modules, we can see that $F(fg)\pi_Z=\pi_Xb_{fg}$ and
$F(fg)\pi_Z=\pi_X(b_fb_g)$. By the uniqueness of
$\underline{b_{fg}}$, we have
$\underline{b_{fg}}=\underline{b_{f}}\;\underline{b_{g}}.$

Moreover, if $X$ is a projective $A$-module, then $F(X)\simeq
\cpx{\bar{Q}_X}$ and $\cpx{\bar{Q}_X}\in \add(\cpx{\bar{Q}})$ by the
proof of Lemma \ref{formofFX}. In particular, $\bar{Q}_X^0$ is
projective. Thus, if $f$ factorizes through a projective module $P$,
say $f=gh$ with $g\in \Hom_A(X,P)$ and $h\in\Hom_A(P,Y)$, then
$b_{f}$ factorizes through a projective $B$-module since $b_f
=b_{gh}=(b_{gh}-b_{g}b_{h})+b_{g}b_{h}$ and since both
$b_{gh}-b_{g}b_{h}$ and  $b_{g}b_{h}$ factorize through projective
$B$-modules.

For each $A$-module $X$, we define $\bar{F}(X)=\bar{Q}_X^0$. Note
that $\bar{Q}_X^0$ is, up to isomorphisms,  uniquely determined by
$X$ (see Lemma \ref{formofFX}). For each homomorphism
$\underline{f}$ in $\StHom_A(X, Y)$, we set
$\bar{F}(\underline{f})=\underline{b_f}$. Then the above discussions
show that $\bar{F}$ is well-defined on Hom-sets and that $\bar{F}$
is a functor from $\stmodcat{A}$ to $\stmodcat{B}$. Note that
$\bar{F}$ is additive since $F$ is additive.

To finish the proof of the lemma, it remains to show that $\pi_X:
F(X)\lra \bar{F}(X)$ is a natural isomorphism in the quotient
category $\Db{B}/\Kb{\pmodcat{B}}$. That the morphism $\pi_X$ is an
isomorphism follows from the fact that $\bar{Q}_X^+$ is isomorphic
to the zero object in $\Db{B}/\Kb{\pmodcat{B}}$. Clearly, $\pi_X$ is
natural in $X$ since we have a commutative diagram
$$\begin{CD}  F(X) @>{\pi_X}>>\bar{Q}_X^0 \\
@V{F(f)}VV @VV{b_{f}}V \\
F(Y)@>{\pi_Y}>>\bar{Q}_Y^0 \end{CD}$$ in the quotient category
$\Db{B}/\Kb{\pmodcat{B}}$. $\dickebox$

\medskip
It is appropriate to introduce a name for the functor $\bar{F}$.
Given a derived equivalence $F$, the functor $\bar{F}$ constructed
in Proposition \ref{stablefunctor} is called a {\em stable functor}
of $F$ throughout this paper.

\begin{Prop}\label{inversestablefunctor}
If $\add(_AQ)=\add(\nu_AQ)$, then there is an additive functor
$\bar{G}: \stmodcat{B}\lra\stmodcat{A}$ sending $U$ to $Q_{U}^0$
such that the following diagram of the additive functors
$$\begin{CD}
\stmodcat{B}@>{\Sigma}>>\Db{B}/\Kb{\pmodcat{B}}\\
@V{\bar{G}}VV   @VV{G}V\\
\stmodcat{A}@>{\Sigma}>>\Db{A}/\Kb{\pmodcat{A}}
\end{CD}$$
is commutative up to natural isomorphisms.
\end{Prop}

{\it Proof.} The idea of the proof of Proposition
\ref{inversestablefunctor} is similar to that of Proposition
\ref{stablefunctor}. We just outline the key points of the
construction of $\bar{G}$.

For a $B$-module $U$, by Lemma \ref{formofGY}, $G(U)$ is isomorphic
in $\Db{A}$ to a complex $\cpx{Q_U}$ such that $Q_U^i\in
\add(\nu_AQ)$ for all $i<0$ and $Q_U^j=0$ for all $j>0$. By the
remark at the end of Section 2, we can assume that $G(U)$ is just
$\cpx{{Q}_U}$. Let $Q_U^-$ denote the complex
$\sigma_{<0}\cpx{{Q}_U}$. We have a distinguished triangle in
$\Db{A}$:
$$\xymatrix{
Q_U^-[-1]\ar[r]^(.60){\beta_U} & Q_U^0\ar[r]^(.4){\lambda_U} &
G(U)\ar[r]^{\gamma_U} & Q_U^-. }$$ Now if $g: U\lra V$ is a
homomorphism of $B$-modules, then we have a commutative diagram
$$\begin{CD}  Q_U^-[-1]@>{\beta_U}>>Q_U^0@>{\lambda_U}>>G(U) @>{\gamma_U}>>Q_U^-\\
@VV{v_g[-1]}V @VV{u_g}V @VV{G(g)}V @VV{v_g}V \\
Q_V^-[-1] @>{\beta_V}>>Q_V^0@>{\lambda_V}>>G(V) @>{\gamma_V}>>Q_V^-.
\end{CD}$$
The existence of $u_g$ follows from the fact that the morphism
$\lambda_UG(g)\gamma_V$ belongs to $\Hom_{\Db{A}}(Q_U^0,Q_V^-)$
$\simeq$ $\Hom_{\Kb{A}}(Q_U^0,Q_V^-)=0$. Since $\modcat{A}$ is fully
embedded in $\D{A}$, the map $u_g$ can be chosen to be an $A$-module
homomorphism. Moreover, if $u_g': Q_U^0\longrightarrow Q_V^0$ is
another morphism such that $u_g'\lambda_V=\lambda_UG(g)$, then
$(u_g-u_g')\lambda_V=0$ and $u_g-u_g'$ factorizes through
$Q_V^-[-1]$. Since $\add(_AQ)=\add(\nu_AQ)$, all the terms of the
complex $Q_V^-[-1]$ are projective-injective. By Lemma
\ref{factthroughP}, the morphism $u_g-u_g'$ factorizes through a
projective module. Thus, for each $B$-module homomorphism $g$, the
morphism $\underline{u_g}$ in $\StHom_A(Q_U^0, Q_V^0)$ is uniquely
determined by $g$.

As in the proof of Proposition \ref{stablefunctor}, we can show that
the composition of two morphisms is preserved, namely
$\underline{u_{gh}}= \underline{u_g}\underline{u_{h}} $ for all
$g\in \Hom_B(U,V)$ and all $h\in \Hom_B(V,W)$.

Moreover, if $P$ is a projective $B$-module, then $\cpx{{Q}_{P}}$ is
isomorphic in $\Db{A}$ to a complex $\cpx{Q_1}$ in $\add(\cpx{Q})$.
Since $\add(_AQ)=\add(\nu_AQ)$, the complex $\cpx{Q_1}$ is of the
form in Lemma \ref{formofGY}. By the uniqueness of $\cpx{Q_P}$, we
have an isomorphism $\cpx{{Q}_{P}}\simeq\cpx{Q_1}$ in $\Cb{A}$.
Hence $Q_P^0\simeq Q_1^0$ and $Q_P^0$ is a projective $A$-module.
Thus, if $g: U\longrightarrow V$ factorizes through a projective
$B$-module $P$, that is, $g=st$ for $s: U\longrightarrow P$ and $t:
P\longrightarrow V$, then
$u_g=u_{st}=(u_{st}-u_{s}u_{t})+u_{s}u_{t}$ factorizes through a
projective $A$-module. This shows that the map $\underline{g}\mapsto
\underline{u_g}$ is well-defined.

\smallskip
Now, we define $\bar{G}(U):= Q_U^0$ for each $B$-module $U$ and
$\bar{G}(\underline{g}):=\underline{u_g}$ for each morphism
$\underline{g}$ in $\stmodcat{B}$. Note that $Q_U^0$ is, up to
isomorphisms, uniquely determined by $U$ (see Lemma \ref{formofGY}).
Thus we obtain an additive functor $\bar{G}$ from $\stmodcat{B}$ to
$\stmodcat{A}$. Moreover, the map $\lambda_U$ is a natural
isomorphism in the quotient category $\Db{A}/\Kb{\pmodcat{A}}$ since
$Q_U^-$ is in $\Kb{\pmodcat{A}}$. $\dickebox$

\begin{Prop}
Suppose $\add (_AQ)=\add(\nu_AQ)$. Let $\bar{F}$ and $\bar{G}$ be
the functors constructed in Proposition \ref{stablefunctor} and
Proposition \ref{inversestablefunctor}, respectively. Then the
composition $\bar{G}\bar{F}$ is naturally isomorphic to the identity
functor $1_{\stmodcat{A}}$. In particular, $\bar{G}$ is dense, and
the restriction of $\bar{G}$ to $\emph{Im}(F)$ is full.
\label{fulldense}
\end{Prop}

{\it Proof.} For each $A$-module $X$, we may assume that $F(X)$ is
the complex $\cpx{\bar{Q}_X}$ defined in Lemma \ref{formofFX}. For
each $B$-module $U$, we assume that $G(U)$ is the complex
$\cpx{Q_U}$ defined in Lemma \ref{formofGY}. We set
$\bar{Q}_X^+=\sigma_{\geq 1}\cpx{\bar{Q}_X}$ and
$Q_U^-=\sigma_{<0}\cpx{{Q}_U}$. Then all the terms of $Q_U^-$ are
projective-injective since $\add(_AQ)=\add(\nu_AQ)$. By definition,
we have $\bar{F}(X)=\bar{Q}_X^0$ for each $A$-module $X$ (see
Proposition \ref{stablefunctor}), and $\bar{G}(U)=Q_U^0$ for each
$B$-module $U$ (see Proposition \ref{inversestablefunctor}). Thus,
for each $A$-module $X$, there is a distinguished triangle
$$\xymatrix{
\bar{Q}_X^+\ar[r]^{i_X} & F(X)\ar[r]^{\pi_X} &
\bar{F}(X)\ar[r]^{\alpha_X} & \bar{Q}_X^+[1] }$$ in $\Db{B}$, and a
distinguished triangle
$$\xymatrix{
Q_{\bar{F}X}^- [-1]\ar[r]^{\beta_{\bar{F}X}} &
\bar{G}\bar{F}(X)\ar[r]^{\lambda_{\bar{F}X}} &
G{\bar{F}(X)}\ar[r]^{\gamma_{\bar{F}X}} & Q_{\bar{F}X}^- }$$ in
$\Db{A}$. Applying $G$ to the first triangle, we obtain the
following commutative diagram in $\Db{A}$
$$\begin{CD}
G\bar{Q}_X^+@>{Gi_X}>>GF(X)@>{G\pi_X}>>G\bar{F}(X)@>{G\alpha_X}>>G\bar{Q}_X^+[1]\\
@VV{q_X[-1]}V @VV{\eta_X}V  @|  @VV{q_X}V\\
Q_{\bar{F}X}^-[-1]@>{\beta_{\bar{F}X}}>>\bar{G}\bar{F}(X)@>{\lambda_{\bar{F}X}}>>G{\bar{F}(X)}
@>{\gamma_{\bar{F}X}}>> Q_{\bar{F}X}^- \,.
\end{CD}$$
The existence of $\eta_X$ follows from the fact that
$G(\pi_X)\gamma_{\bar{F}X}$ belongs to $\Hom_{\Db{A}}(GF(X),
{Q_{\bar{F}X}^-})\simeq \Hom_{\Db{A}}(X, {Q_{\bar{F}X}^-})$ $\simeq
\Hom_{\Kb{A}}(X, {Q_{\bar{F}X}^-})=0$. Since $GF$ is naturally
isomorphic to the identity functor $1_{\Db{A}}$, there is a natural
morphism $\epsilon_X: X\lra GF(X)$ in $\Db{A}$ for each $A$-module
$X$. Let $\theta_X$ be the composition $\epsilon_X\eta_X$. Then
$\theta_X: X\lra \bar{G}\bar{F}(X)$ is an $A$-module homomorphism
since $\modcat{A}$ is fully embedded in $\Db{A}$.

We claim that $\underline{\theta_X}$ is a natural map in
$\stmodcat{A}$. Indeed, for any $A$-module homomorphism $f:
X\rightarrow Y$, by the proof of Proposition \ref{stablefunctor}, we
have a homomorphism $b_f: \bar{F}(X)\lra\bar{F}(Y)$ of $B$-modules
such that $\pi_Xb_f=F(f)\pi_Y$ in $\Db{B}$. By the proof of
Proposition \ref{inversestablefunctor}, there is a homomorphism
$u_{b_f}: \bar{G}(\bar{F}(X))\lra \bar{G}(\bar{F}(Y))$ of
$A$-modules such that
$u_{b_f}\lambda_{\bar{F}Y}=\lambda_{\bar{F}X}G(b_f)$ in $\Db{A}$.
Thus, we have in $\Db{A}$:

$$\begin{array}{rl}
(\theta_Xu_{b_f}-f\theta_Y)\lambda_{\bar{F}Y} \medskip & =
(\epsilon_X\eta_Xu_{b_f}-f\epsilon_Y\eta_Y) \lambda_{\bar{F}Y}\\
\medskip
&=(\epsilon_X\eta_Xu_{b_f}-\epsilon_XGF(f)\eta_Y)\lambda_{\bar{F}Y}\\
\medskip
&
=\epsilon_X(\eta_Xu_{b_f}\lambda_{\bar{F}Y}-GF(f)\eta_Y\lambda_{\bar{F}Y})\\
\medskip
&
=\epsilon_X(\eta_X\lambda_{\bar{F}X}G(b_f)-GF(f)\eta_Y\lambda_{\bar{F}Y})\\
\medskip
& =\epsilon_X(G(\pi_X)G(b_f)-GF(f)G(\pi_Y))\\ \medskip &
=\epsilon_X(G(\pi_X)G(b_f)-G(F(f)\pi_Y))\\ \medskip
& =\epsilon_X(G(\pi_X)G(b_f)-G(\pi_Xb_f))\\ \medskip  & =0.\\
\end{array}$$
This implies that the map $\theta_Xu_{b_f}-f\theta_Y$ factorizes
through ${Q_{\bar{F}Y}^-}[-1]$. It follows by Lemma
\ref{factthroughP} that $\theta_Xu_{b_f}-f\theta_Y$ factorizes
through a projective module. Note that
$\underline{u_{b_f}}=\bar{F}\bar{G}(\underline{f})$. Thus
$\underline{\theta_X}\bar{F}\bar{G}(\underline{f})-\underline{f}\underline{\theta_Y}=0$
in $\stmodcat{A}$ and $\underline{\theta_X}$ is natural in $X$.

To finish the proof, we have to show that $\underline{\theta_X}$ is
an isomorphism in $\stmodcat{A}$ for each $A$-module $X$. Clearly,
we can assume that $X$ is an indecomposable non-projective
$A$-module. Using the method similar to that in the proof of Lemma
\ref{formofFX}, we can prove that $G(\bar{Q}_X^+)$ is isomorphic in
$\Db{A}$ to a radical complex $\cpx{Q_1}$ in $\Kb{\pmodcat{A}}$ with
$Q_1^i\in\add(_AQ)$ for all $i\leq 0$. Since both $X$ and
$G\bar{F}(X)$ have no homology in positive degrees, the complex
$G({\bar{Q}_X^+})$ has no homology in degrees greater than $1$, and
therefore $Q_1^i=0$ for all $i>1$. Now we may form the following
commutative diagram in $\Db{A}$:
$$\begin{CD}
\cpx{Q_1}@>{\phi_X}>> X@>{\lambda}>>\mbox{con}(\phi_X)@>>>\cpx{Q}_1[1]\\
@VV{s}V  @VV{\epsilon_X}V  @VV{t}V  @VV{s[1]}V\\
G({\bar{Q}_X^+})@>{Gi_X}>>GF(X)@>{G\pi_X}>>G\bar{F}(X)@>>> G({\bar{Q}_X^+})[1]\\
@.    @VV{\eta_X}V @| \\
@.  \bar{G}\bar{F}X@>{\lambda_{\bar{F}X}}>>G\bar{F}(X),
\end{CD}$$
where $s$ is an isomorphism from $G(Q_X^+)$ to $\cpx{Q_1}$, and
where $\phi_X=sG(i_X)\epsilon^{-1}$ is a chain map, and where
$\lambda$ is induced by the canonical map $\lambda^0$ from $X$ to
$Q_1^1\oplus X$ defined by the mapping cone. Since $Q_1^i$ is in
$\add(_AQ)$ for all $i\leq 0$ and zero for all $i>1$, the mapping
cone $\mbox{con}(\phi_X)$ has terms in $\add(_AQ)$ for all negative
degrees and zero for all positive degrees. Note that
$\add(_AQ)=\add(\nu_AQ)$ and the $A$-module $Q$ is
projective-injective. Consequently, the  terms of
$\mbox{con}(\phi_X)$ in all negative degrees are
projective-injective. Note that
$G(\bar{F}(X))\simeq\cpx{Q_{\bar{F}(X)}}$ by our assumption. Thus,
by definition (see Lemma \ref{formofGY}), all terms of
$\cpx{Q_{\bar{F}(X)}}$ in negative degrees are projective-injective.
Moreover, both $\mbox{con}(\phi_X)$ and $\cpx{Q_{\bar{F}(X)}}$ have
zero terms in all positive degrees. By Lemma \ref{kdiso} (2), we
have $\Hom_{\Db{A}}(\mbox{con}(\phi_X),\cpx{Q_{\bar{F}(X)}})=
\Hom_{\Kb{A}}(\mbox{con}(\phi_X),\cpx{Q_{\bar{F}(X)}})$.  Since the
two morphisms $\epsilon_X$ and $s$ are both isomorphism in $\Db{A}$,
the morphism $t$ is also an isomorphism in $\Db{A}$. Hence $t$ is an
isomorphism from $\mbox{con}(\phi_X)$ to $\cpx{Q_{\bar{F}(X)}}$ in
$\Kb{A}$. Moreover, since $X$ is indecomposable and non-projective,
the complex $\mbox{con}(\phi_X)$ is a radical complex. Thus, the
chain map $t$ is actually an isomorphism between
$\mbox{con}(\phi_X)$ and $\cpx{Q_{\bar{F}X}}$ in $\Cb{A}$, and the
morphism $t^0: (\mbox{con}(\phi_X))^0=Q_1^1\oplus X\longrightarrow
Q_{\bar{F}X}^0$ in  degree zero is an isomorphism of $A$-modules.
From the above commutative diagram, we have
$\theta_X\lambda_{\bar{F}X}-\lambda t=0$ in $\Db{A}$. By Lemma
\ref{kdiso} (2), we see that $\theta_X\lambda_{\bar{F}X}-\lambda t$
is null-homotopic in $\Cb{A}$. This means that $\theta_X-\lambda^0
t^0$ factorizes through the projective $A$-module
$Q_{\bar{F}X}^{-1}$. Hence
$\underline{\theta_X}=\underline{\lambda}^0\underline{t}^0$ is an
isomorphism in $\stmodcat{A}$ since $\underline{\lambda^0}$ and
$\underline{t^0}$ both are isomorphisms in $A$-{\underline{mod}.
$\dickebox$

\medskip
{\it Remark.} Without the condition $\add(_AQ)=\add(\nu_AQ)$ in
Proposition \ref{inversestablefunctor}, we can similarly define a
functor $\bar{G}': \stmodcat{B}\lra\stmodcat{A}$, as was done in
Proposition \ref{stablefunctor}. But the disadvantage of using
$\bar{G}'$ is that we do not know any behavior of the composition of
$\bar{F}$ with $\bar{G}'$.

\medskip
\medskip
We say that a derived equivalence $F$ between Artin algebras $A$ and
$B$ is {\it almost $\nu$-stable} if $\add(_AQ)=\add(\nu_AQ)$ and
$\add(_B\bar{Q})=\add(\nu_B\bar{Q}).$
\medskip

The following theorem shows that the almost $\nu$-stable condition
is sufficient for $\bar{F}$ to be an equivalence.

\begin{Theo}
Let $A$ and $B$ be two Artin $R$-algebras, and let
$F:\Db{A}\lra\Db{B}$ be a derived equivalence. If $F$ is almost
$\nu$-stable, then the stable functor $\bar{F}$ is an equivalence.
\label{Derivedstable}
\end{Theo}

{\it Proof.} Since $F$ is almost $\nu$-stable, we have
$\add(_AQ)=\add(\nu_AQ)$. By Proposition \ref{fulldense}, we have
$\bar{G}\bar{F}\simeq 1_{\stmodcat{A}}$. Since $F$ is almost
$\nu$-stable, we also have $\add(\bar{Q})=\add(\nu\bar{Q})$. With a
proof similar to that of Proposition \ref{fulldense}, we can show
that $\bar{F}\bar{G}$ is naturally isomorphic to the identity
functor $1_{\stmodcat{B}}$. Thus $\bar{F}$ and $\bar{G}$ are
equivalences of categories. $\dickebox$

\medskip
Theorem \ref{Derivedstable} gives rise to a method of getting stable
equivalences from derived equivalences. In Section \ref{sectStM}, we
shall prove that, for finite-dimensional algebras, one even can get
a stable equivalence of Morita type, which has many pleasant
properties (see \cite{Broue1994}, \cite{Xi1}, \cite{Xi2} and the
references therein).

\medskip
In the following, we shall develop some properties of almost
$\nu$-stable functors, which will be used in Section \ref{sectStM}.

Let $_AE$ be a direct sum of all those non-isomorphic indecomposable
projective-injective $A$-modules $X$ that have the property:
$\nu_A^iX$ is again projective-injective for every $i> 0$. The
$A$-module $_AE$ is unique up to isomorphism, and it is called the
{\em maximal $\nu$-stable} $A$-module. Similarly, we have a maximal
$\nu$-stable $B$-module $_B\bar{E}$. The following result shows that
an almost $\nu$-stable functor is closely related to the maximal
$\nu$-stable modules.

\begin{Prop} The following are equivalent:

  $(1)$ $F$ is almost $\nu$-stable, that is, $\add(\nu_AQ)=\add(_AQ)$ and
$\add(\nu_B\bar{Q})=\add(_B\bar{Q})$.

  $(2)$ $_AQ\in\add(_AE)$ and $_B\bar{Q}\in\add(_B\bar{E})$.

  $(3)$ $_AQ$ and $\nu_B\bar{Q}$ are projective-injective.  \label{prop4.1}
\end{Prop}
{\it Proof. } Clearly, we have $(1)\Rightarrow(2)\Rightarrow(3)$.
Now we show that $(3)$ implies $(1)$. Assume that $_AQ$ is
injective. By Lemma \ref{formofGY}, $G(B)$ is isomorphic in $\Db{A}$
to a radical complex $\cpx{Q_B} =(Q_B^i, d^i)$ with $Q_B^i$ in
$\add(\nu_AQ)$ for all $i<0$. In particular, $Q^i_B$ is
projective-injective for all $i<0$. Since $G(B)\simeq\cpx{Q}$, the
complexes $\cpx{Q}$ and $\cpx{Q_B}$ are isomorphic in $\Db{A}$.
Since $_AQ$ is injective by assumption, all the terms of $\cpx{Q}$
in negative degrees are injective. By Lemma \ref{kdiso}, the
complexes $\cpx{Q}$ and $\cpx{Q_B}$ is isomorphic in $\Kb{A}$. Since
both $\cpx{Q}$ and $\cpx{Q_B}$ are radical, they are also isomorphic
in $\Cb{A}$. In particular, we have $Q^i\simeq Q_B^i$ for all $i<0$,
and therefore $_AQ
:=\bigoplus_{i=-1}^{-n}Q^i\simeq\bigoplus_{i=-1}^{-n}Q_B^i\in\add(\nu_AQ)$.
Since $_AQ$ and $\nu_AQ$ has the same number of isomorphism classes
of indecomposable direct summands, we have $\add(_AQ)=\add(\nu_AQ)$.
Similarly, we have $\add(_B\bar{Q})=\add(\nu_B\bar{Q})$. This proves
$(3)\Rightarrow(1)$.
 $\dickebox$

\medskip
{\it Remark: } From Lemma \ref{prop4.1}, we can see that every
derived equivalence between two self-injective Artin algebras is
almost $\nu$-stable. Thus, we can re-obtain the result
\cite[Corollary 2.2]{RickDstable} of Rickard by Theorem
\ref{Derivedstable}.

\begin{Lem}
Suppose $Q\in \add(_AE)$, and $\bar{Q}\in \add(_B\bar{E})$. Then

$(1)$  for each $\cpx{P}$ in $\Kb{\add(_AE)}$, the complex
$F(\cpx{P})$ is isomorphic in $\Db{B}$ to a complex in
$\Kb{\add(_B\bar{E})}$.

$(2)$ for each $\cpx{\bar{P}}$ in $\Kb{\add(_B\bar{E})}$, the
complex $G(\cpx{\bar{P}})$ is isomorphic in $\Db{A}$ to a complex in
$\Kb{\add(_AE)}$. \label{lem4.2}
\end{Lem}
{\it Proof.} (1) It suffices to show that, for each indecomposable
$A$-module $U$ in $\add(_AE)$, the complex $F(U)$ is isomorphic to a
complex in $\Kb{\add(_B\bar{E})}$. Suppose $U\in\add(_AE)$. By Lemma
\ref{formofFX}, $F(U)$ is isomorphic in $\Db{B}$ to a radical
complex $\cpx{\bar{Q}_U}$:
$$\xymatrix{
0 \ar[r] & \bar{Q}_U^0 \ar[r] & \bar{Q}_{U}^1\ar[r]& \cdots \ar[r] &
\bar{Q}_{U}^n\ar[r]& 0 }$$ with $\bar{Q}_U^i\in\add(_B\bar{Q})$ for
all $i>0$. For simplicity, we assume that $F(U)$ is just
$\cpx{\bar{Q}_U}$. Since
$\bar{Q}_U^i\in\add(_B\bar{Q})\subseteq\add(_B\bar{E})$ for $i>0$,
it remains to show that $\bar{Q}_U^0$ is in $\add(_B\bar{E})$.
Clearly, $\bar{Q}_U^0$ is projective since $U$ is projective. Note
that we have an isomorphism $\nu_BF(U)\simeq F(\nu_AU)$ in $\Db{B}$,
that is, $\nu_B\cpx{\bar{Q}_U}$ is isomorphic to
$\cpx{\bar{Q}_{\nu_AU}}$. Note that $\nu_B\bar{Q}_U^i\in
\add(\nu_B\bar{Q})$ for all $i>0$ and $\add(\nu_B\bar{Q})\subseteq
\add(_B\bar{E})$ by the definition of $_B\bar{E}$. Thus
$\nu_B\bar{Q}_U^i$ is projective-injective for all $i>0$, and
$\nu_B\cpx{\bar{Q}_U}$ and $\cpx{\bar{Q}_{\nu_AU}}$ are isomorphic
in $\Kb{B}$ by Lemma \ref{kdiso}. Since both $\nu_B\cpx{\bar{Q}_U}$
and $\cpx{\bar{Q}_{\nu_AU}}$ are radical complexes,
$\nu_B\cpx{\bar{Q}_U}$ and $\cpx{\bar{Q}_{\nu_AU}}$ are actually
isomorphic in $\Cb{B}$, and particularly we have
$\nu_B\bar{Q}_U^0\simeq\bar{Q}_{\nu_AU}^0$. Note that if
$U\in\add(_AE)$, then $\nu^i_AU\in \add(E)$ for all $i\geq 0$ by
definition. Hence, for each integer $m>0$, we have
$\nu_B^m(\bar{Q}_U^0)\simeq\nu^{m-1}_B(\bar{Q}_{\nu_AU}^0)\simeq\cdots\simeq\bar{Q}_{\nu_A^mU}^0$,
and therefore $\nu_B^m(\bar{Q}_U^0)$ is projective-injective. Thus
 $\bar{Q}_U^0\in\add(_B\bar{E})$ by definition.

(2) is a dual statement of (1).  $\dickebox$

\medskip
The following is a consequence of Lemma \ref{lem4.2}.

\begin{Koro} If $F$ is an almost $\nu$-stable derived equivalence between Artin
algebras $A$ and $B$, then there is a derived equivalence between
the self-injective algebras $\End_A(E)$ and $\End_B(\bar{E})$.
\label{3.10}
\end{Koro}

{\it Proof.} By Lemma \ref{lem4.2}, $F$ induces an equivalence
between $\Kb{\add(_AE)}$ and $\Kb{\add(_B\bar{E})}$ as triangulated
categories. Since $\Kb{\add(_AE)}$ and $\Kb{\pmodcat{\End(_AE)}}$
are equivalent as triangulated categories, we obtain an equivalence
between $\Kb{\pmodcat{\End(_AE)}}$ and
$\Kb{\pmodcat{\End(_B\bar{E})}}$ as triangulated categories. By
\cite[Theorem 6.4]{RickMoritaTh}, the algebras $\End_A(E)$ and
$\End_B(\bar{E})$ are derived-equivalent. Note that $\End_A(E)$ is
self-injective.
$\dickebox$

Let us end this section by the following result which tells us how
to get an almost $\nu$-stable derived equivalence from a tilting
module.

Let $A$ be an Artin algebra. Recall that an $A$-module $T$ is called
a \emph{tilting module} if

$(1)$ the projective dimension of $T$ is finite,

$(2)$ $\Ext^i_A(T,T)=0$ for all $i>0$, and

$(3)$ there is an exact sequence $0\lra A\lra T^0\lra\cdots\lra
T^m\lra 0$ of $A$-mod with each $T^i$ in $\add(_AT)$.

\medskip
It is well-known that a tilting $A$-module $_AT$ induces a derived
equivalence between $A$ and $\End_A(T)$ (see \cite[Theorem 2.10, p.
109]{HappelTriangle} and \cite[Theorem 2.1]{ClineParshallScott}).

\begin{Prop}
Let $A$ be an Artin algebra. Suppose $_AT$ is a tilting $A$-module
with $B=\End_A(T)$. Let
$$0\lra P_n\stackrel{d_n}{\lra} P_{n-1}\lra\cdots\lra P_0\stackrel{d_0}{\lra} T\lra 0$$
be a minimal projective resolution of $_AT$. Set
$P:=\bigoplus_{i=0}^{n-1}P_i$. If $\add(_AP)=\add(\nu_AP)$, then
there is an almost $\nu$-stable derived equivalence between $A$ and
$B$. \label{stabletilting}
\end{Prop}

{\it Proof: } By \cite[Theorem 2.1]{ClineParshallScott}, the functor
$F'={}_AT\otimesL_B-: \Db{B}\lra\Db{A}$ is a derived equivalence.
Now we denote $F'[-n]$ by $F$. Let $\cpx{P}$ be the complex
$$0\lra P_n\stackrel{d_n}{\lra} P_{n-1}\lra\cdots\lra P_0\lra 0$$
with $P_n$ in degree zero. Then we have
$F(B)=(_AT\otimesL_BB)[-n]={}_AT[-n]\simeq \cpx{P}$ in $\Db{A}$. Let
$G$ be a quasi-inverse of $F$. Then $G(\cpx{P})\simeq G(F(B))\simeq
B$ in $\Db{B}$, and therefore $\cpx{P}$ is a radical tilting complex
associated to $G$.

Since $\add(\nu_AP)=\add(_AP)$, the module $_AP$ is
projective-injective. Thus $_AP\in\add(_AT)$ and $P_i\in\add(_AT)$
for all $0\le i\le n-1$. We denote by $\cpx{T}$ the complex
$$0\lra P_{n-1}\oplus P\stackrel{\left[{{d_{n-1}}\atop{0}}\right]}{\lra} P_{n-2}
\lra\cdots\lra P_0\stackrel{d_0}{\lra} T\lra 0$$ with $T$ in degree
zero. Then $\HomP_A(T, \cpx{T})$ is a complex in $\Kb{\pmodcat{B}}$,
and
$$F(\HomP_A(T, \cpx{T}))={}_AT\otimesL_B\HomP_A(T, \cpx{T})[-n]\simeq
{}_AT\otimesP_B\HomP_A(T, \cpx{T})[-n]\simeq \cpx{T}[-n]\simeq
P_n\oplus P$$
 in $\Db{A}$. Since $\cpx{P}$ is a tilting complex over $A$, we have
$_AA\in\add(P_n\oplus P)$. Thus there is a radical complex
$\cpx{\bar{P}}$ in $\Kb{\pmodcat{B}}$ such that
$\cpx{\bar{P}}\in\add(\HomP_A(T, \cpx{T}))$ and
$F(\cpx{\bar{P}})\simeq A$ in $\Db{A}$. By definition,
$\cpx{\bar{P}}$ is a tilting complex associated to $F$. (For the
unexplained notations appearing in this proof, we refer the reader
to Section \ref{sectStM} below).

We claim that $F$ is almost $\nu$-stable. In fact,
$\bigoplus_{i=-1}^{-n}\bar{P}^i$ is in $\add(\Hom_A(T, P))$, and
$\bigoplus_{i=1}^{n}P^i=\bigoplus_{i=0}^{n-1}P_i=P$. Let $_AE$
(respectively, $_B\bar{E}$) be the maximal $\nu$-stable $A$-module
(respectively, $B$-module). Then $_AP\in\add(_AE)$. Note that we
have the following isomorphisms of $B$-modules:
$$\nu_B\Hom_A(T, P)=D\Hom_B(\Hom_A(T, P), \Hom_A(T, T))\simeq D\Hom_A(P, T)
\simeq D(P^*\otimes_AT)\simeq\Hom_A(T, \nu_AP).$$ Since
$\add(_AP)=\add(\nu_AP)$, we have $\add(\nu_B\Hom_A(T,
P))=\add(\Hom_A(T, P))$, that is, $\Hom_A(T, P)\in\add(_B\bar{E})$.
It follows that $\bigoplus_{i=-1}^{-n}\bar{P}^i$ is in
$\add(_B\bar{E})$. By Proposition \ref{prop4.1}, the functor $F$ is
almost $\nu$-stable.
 $\dickebox$

 \medskip
{\it Remark: } Let $A$ be a self-injective Artin algebra, and let
$X$ be an $A$-module. By \cite[Corollary 3.7]{HuXi2}, there is a
derived equivalence between $\End_A(A\oplus X)$ and
$\End_A(A\oplus\Omega_A(X))$ induced by the almost $\add(A)$-split
sequence $0\ra \Omega_A(X)\ra P_X\ra X\ra 0$, where $P_X$ is a
projective cover of $X$. By Proposition \ref{stabletilting}, it is
easy to check  that this is an almost $\nu$-stable derived
equivalence. Thus the algebras $\End_A(A\oplus X)$ and
$\End_A(A\oplus\Omega_A(X))$ are stably equivalent by Theorem
\ref{Derivedstable}.

\section{Comparison of homological dimensions \label{dimensions}}
In this section, we shall deduce some basic homological properties
of the functor $\bar{F}$ and compare homological dimensions of $A$
with that of $B$.

We make the following convention: From now on, throughout this
paper, we keep our notations introduced in the previous sections.

Recall that the \emph{finitistic dimension} of an Artin algebra $A$,
denoted by fin.dim$(A)$, is defined to be the supremum of the
projective dimensions of finitely generated $A$-modules of finite
projective dimension. The finitistic dimension conjecture states
that fin.dim$(A)$ should be finite for any Artin algebra $A$.
Concerning the new advances on this conjecture, we refer the reader
to the recent paper \cite{xx} and the references therein.

For an $A$-module $X$, we denote by $\pd(_AX)$ the projective
dimension of $X$, and by gl.dim$(A)$ the global dimension of $A$,
which is by definition the supremum of the projective dimensions of
all finitely generated $A$-modules. Clearly, if gl.dim$(A)<\infty$,
then fin.dim$(A)$ = gl.dim$(A)$.

\begin{Prop}
Let $\bar{F}$ be the stable functor of $F$ defined in Proposition
\ref{stablefunctor}. Then:

$(1)$ For each exact sequence $0\longrightarrow
X\stackrel{f}{\longrightarrow} Y\stackrel{g}{\longrightarrow}
Z\longrightarrow 0$ in $\modcat{A}$, there is an exact sequence
$$\begin{CD}
0@>>>\bar{F}(X)@>{[a,f']}>>P\oplus\bar{F}(Y)@>{\left[{{b}\atop{g'}}\right]}>>\bar{F}(Z)@>>>0
\end{CD}$$ in $B\modcat$, where $P\in\add(_B\bar{Q})$,
$\bar{F}(\underline{f})=\underline{f}'$ and
$\bar{F}(\underline{g})=\underline{g}'$.

$(2)$ For each $A$-module $X$, we have a $B$-module isomorphism:
$\bar{F}(\Omega_A(X))\simeq \Omega_B(\bar{F}(X))\oplus P$, where $P$
is a projective $B$-module, and $\Omega$ is the syzygy operator.

$(3)$ For each $A$-module $X$, we have $\pd(_B\bar{F}(X))\leq \pd(_A
X)\leq \pd(_B\bar{F}(X))+n$.

$(4)$ If $\bar{F}$ is an equivalence, then $A$ and $B$ have the same
finitistic and global dimensions.\label{barFProp}
\end{Prop}

{\it Proof.} For each $A$-module $X$, we may assume that $F(X)$ is
the complex $\cpx{\bar{Q}_X}$ defined in Lemma \ref{formofFX}.

(1) From the exact sequence $0\longrightarrow
X\stackrel{f}{\longrightarrow} Y\stackrel{g}{\longrightarrow}
Z\longrightarrow 0$ in $A$-mod we have a distinguished triangle in
$\Db{A}$:
 $$\begin{CD}
    X@>f>> Y@>g>> Z@>\epsilon>> X[1]
 \end{CD}.$$
Applying the functor $F$, we get a distinguished triangle
 $$\begin{CD}
    \cpx{\bar{Q}_X}@>F(f)>> \cpx{\bar{Q}_Y}@>F(g)>> \cpx{\bar{Q}_Z}@>F(\epsilon)>> \cpx{\bar{Q}_X}[1]
 \end{CD}$$
in $\Db{B}$. Moreover, by Lemma \ref{kdiso}, the morphisms $F(f)$
and $F(g)$ are induced by chain maps $\cpx{p}$ and $\cpx{q}$,
respectively. So, we may assume that $F(f)=\cpx{p}$ and
$F(g)=\cpx{q}$. Let $\con(\cpx{q})$ be the mapping cone of the chain
map $\cpx{q}$. Then we have a commutative diagram in $\Db{B}$
$$\begin{CD}
\cpx{\bar{Q}}_Z[-1]@>>>\cpx{\bar{Q}}_X@>{\cpx{p}}>>\cpx{\bar{Q}}_Y@>{\cpx{q}}>> \cpx{\bar{Q}}_Z\\
@|  @VsVV @| @|\\
\cpx{\bar{Q}}_Z[-1]@>>>\con(\cpx{q})[-1]@>{\cpx{\pi}}>>\cpx{\bar{Q}}_Y@>{\cpx{q}}>> \cpx{\bar{Q}}_Z\\
\end{CD}$$
for some isomorphism $s$, where $\cpx{\pi}=(\pi^i)$ with $\pi^i:
\bar{Q}_Y^i\oplus \bar{Q}_Y^{i-1}\longrightarrow \bar{Q}_Y^i$  the
canonical projection for each integer $i$. By Lemma \ref{kdiso}, the
morphism $s$ is induced by a chain map $\cpx{s}$.  By definition,
$\con(\cpx{s})$ is the following complex
$$\begin{CD}
0@>>> \bar{Q}_X^0@>{\left[-d, s^0\right]}>> \bar{Q}_X^1\oplus
\bar{Q}_Y^0@>{\left[\begin{array}{ccc} -d & u & v\\
0 & -d & q^0\end{array}\right]}>> \bar{Q}_X^2\oplus
\bar{Q}_Y^1\oplus\bar{Q}_Z^0@>>>\cdots @>>>
 \bar{Q}_Z^{n}@>>> 0,
\end{CD} $$
where $s^1=[u, v]: \bar{Q}_X^1\longrightarrow
\bar{Q}_Y^1\oplus\bar{Q}_Z^0$, and where the modules $\bar{Q}_X^i$,
$\bar{Q}_Y^i$ and $\bar{Q}_Z^i$ are projective for $i>0$. Since $s$
is an isomorphism in $\Db{B}$, the mapping cone $\con(\cpx{s})$ of
$\cpx{s}$ is an  acyclic complex. Note that the map $-d:
\bar{Q}_Y^0\longrightarrow \bar{Q}_Y^1$ is a radical map. Thus,
dropping the split direct summands of the acyclic complex
$\con(\cpx{s})$, we get an exact sequence
$$\begin{CD}
(*)\qquad 0@>>>\bar{Q}_X^0@>{[a,\; s^0]}>>P\oplus\bar{Q}_Y^0@>{[
b,\; q^0]^T}>>\bar{Q}_Z^0@>>>0
\end{CD}$$
in $B$-mod, where $P$ is a direct summand  of $\bar{Q}_X^1$, and
where $a$ and $b$ are some homomorphisms of $B$-modules. It follows
from $\cpx{p}-\cpx{s}\cpx{\pi}=0$ that the morphism
$\cpx{p}-\cpx{s}\cpx{\pi}$ is null-homotopic by Lemma \ref{kdiso}.
Therefore $p^0-s^0$ factorizes through $\bar{Q}_X^1$, and
$\underline{p}^0=\underline{s}^0$. By setting $f'=s^0$ and $g'=q^0$,
we re-write ($*$) as
$$\begin{CD}
0@>>>\bar{F}(X)@>{[a,f']}>>P\oplus\bar{F}(Y)@>{\left[{{b}\atop{g'}}\right]}>>\bar{F}(Z)@>>>0
\end{CD}$$
with $P\in\add(_B\bar{Q})$, $\bar{F}(\underline{f})=\underline{f}'$
and $\bar{F}(\underline{g})=\underline{g}'$. This proves $(1)$.

\smallskip
$(2)$ Let $X$ be an $A$-module. We have an exact sequence
$0\lra\Omega_A(X)\lra P_X\lra X\rightarrow 0$ in $A$-mod with $P_X$
projective. By $(1)$, we get an exact sequence
$0\lra\bar{F}(\Omega_A(X))\longrightarrow
P\oplus\bar{F}(P_X)\longrightarrow \bar{F}(X)\lra 0$ in $B$-mod for
some projective $B$-module $P$. By the definition of $\bar{F}$, the
$B$-module $\bar{F}(P_X)$ is projective. Thus (2) follows.

$(3)$ The inequality $\pd(_B\bar{F}(X))\leq \pd(_A X)$ follows from
(2). In fact, we may assume $\pd(_AX) = m<\infty$. Then
$\Omega^m_A(X)$ is projective. Therefore $\Omega_B^m\bar{F}(X)$ is
projective by (2), and $\pd(_B\bar{F}(X))\leq m$.

For the second inequality in (3), we may assume $\pd(_B\bar{F}(X))=
m <\infty$. Let $Y$ be an $A$-module. We claim that
$$\Ext_A^i(X, Y)\simeq\Hom_{\Db{A}}(X,Y[i])\simeq \Hom_{\Db{B}}(F(X),F(Y)[i])=0$$
for all $i> m + n$. Indeed, by Lemma \ref{formofFX}, the complex
$F(X)$ is isomorphic in $\Db{B}$ to a complex $\cpx{\bar{Q}_X}$ with
$Q_X^i$ being projective for all $i>0$.  Since
$\pd(_B\bar{Q}_X^0)=\pd(_B\bar{F}(X))=m$, we see that
$\cpx{\bar{Q}_X}$ is isomorphic in $\Db{B}$ to a complex $\cpx{P}$
in $\Kb{\pmodcat{B}}$ with $P^k=0$ for all $k < -m$. Note that
$F(Y)$ is isomorphic to the complex $\cpx{\bar{Q}_Y}$ with
$\bar{Q}_Y^k=0$ for all $k>n$. Clearly, we have
$\Hom_{\Db{B}}(F(X),F(Y)[i])\simeq
\Hom_{\Db{B}}(\cpx{P},\cpx{\bar{Q}_Y}[i])=\Hom_{\Kb{B}}(\cpx{P},\cpx{\bar{Q}_Y}[i])=0$
for all $i>m+n$. Then the second inequality follows.

$(4)$ is a consequence of $(2)$. In fact, suppose $\bar{F}$ is an
equivalence. Then, for an $A$-module $X$ and a positive integer $m$,
the $A$-module $\Omega_A^m(X)$ is projective if and only if
$\bar{F}(\Omega_A^m(X))$ is projective. By (2),
$\bar{F}(\Omega_A^m(X))$ is projective if and only if
$\Omega_B^m(\bar{F}(X))$ is projective. It follows that
$\pd(_AX)\leq m$ if and only if $\pd(_B\bar{F}(X))\leq m$, and
consequently $\pd(_AX)=\pd(_B\bar{F}(X))$ for arbitrary $A$-module
$X$. Thus $(4)$ follows. $\dickebox$

\medskip
{\it Remark.} Proposition \ref{barFProp} (3) can be regarded as an
alternative proof of the fact that if two Artin algebras $A$ and $B$
are derived-equivalent then fin.dim$(A)<\infty$ if and only if
fin.dim$(B)<\infty$.

In Proposition \ref{inversestablefunctor}, we have constructed a
functor $\bar{G}: \stmodcat{B}\longrightarrow\stmodcat{A}$ under the
condition $\add(_AQ)=\add(\nu_AQ)$. This functor $\bar{G}$ has many
properties similar to that of $\bar{F}$.

\begin{Prop}
Suppose $\add(_AQ)=\add(\nu_AQ)$. Let $\bar{G}:
\stmodcat{B}\longrightarrow\stmodcat{A}$ be the functor defined in
Proposition \ref{inversestablefunctor}. Then:

$(1)$ For each exact sequence $0\longrightarrow
U\stackrel{f}{\longrightarrow} V\stackrel{g}{\longrightarrow}
W\longrightarrow 0$ in $\modcat{B}$, there is an exact sequence
$$\begin{CD}
0@>>>\bar{G}(U)@>{[a,f']}>>P\oplus\bar{G}(V)@>{\left[{{b}\atop{g'}}\right]}>>\bar{G}(W)@>>>0
\end{CD}$$
in $A\modcat$, where $P\in\add(_AQ)$,
$\bar{G}(\underline{f})=\underline{f}'$ and
$\bar{G}(\underline{g})=\underline{g}'$.

\smallskip
$(2)$ For each $B$-module $Y$, we have $\bar{G}(\Omega_B(Y))\simeq
\Omega_A(\bar{G}(Y))\oplus P$ in $\modcat{A}$ for a projective
$A$-module $P$.

$(3)$ For each $B$-module $Y$, we have $\pd(_A\bar{G}(Y))\leq \pd(_B
Y)$.

$(4)$ If $I$ is an injective $B$-module, then $\bar{G}(I)$ is an
injective $A$-module. Moreover, $\bar{G}(D(B))\simeq \nu_AQ^0$.
\label{barGProp}
\end{Prop}
{\it Proof: } $(1)$, $(2)$ and $(3)$ are dual statements of
Proposition \ref{barFProp}, and their proofs will be omitted here.
We only prove $(4)$. Let $I$ be an injective $B$-module. Then
$I=\nu_BP$ for a projective $B$-module $P$. Since
$G(B)\simeq\cpx{Q}$, we know that $G(P)$ is isomorphic in
$\Kb{\pmodcat{A}}$ to a radical complex $\cpx{Q_1}\in\add(\cpx{Q})$.
Thus the complex $\cpx{Q_I}$ defined in Lemma \ref{formofGY} is
isomorphic to $G(I)\simeq \nu_AG(P)\simeq\nu_A\cpx{Q_1}$ by Lemma
\ref{nakayama}. Moreover, all the terms of $\nu_A\cpx{Q_1}$ in
negative degrees are in $\add(\nu_AQ)$. Thus, by Lemma
\ref{formofGY}, the complexes $\cpx{Q_I}$ and $\nu_A\cpx{Q_1}$  are
isomorphic in $\Cb{A}$, and consequently $\bar{G}(I)=Q_I^0\simeq
\nu_AQ_1^0$ is injective. In particular, if we take $P={}_BB$ and
$\cpx{Q_1} = \cpx{Q}$, then
$\bar{G}(D(B))=\bar{G}(\nu_BB)\simeq\nu_AQ^0$. $\dickebox$

Let $X$ be an $A$-module, and let $0\lra X\lra I_0\lra I_1\lra\cdots
$ be a minimal injective resolution of $X$ with all $I_j$ injective.
The {\em dominant dimension} of $X$, denoted by $\domdim(X)$, is
defined to be
$$\domdim(X):=\sup\{m \mid I_i \mbox{ is projective for all } 0\le i\leq m-1\}.$$
The \emph{dominant dimension} of the algebra $A$, denoted by
$\domdim(A)$, is defined to be the dominant dimension of the module
$_AA$. Concerning the dominant dimension of an Artin algebra, there
is a conjecture, namely the Nakayama conjecture, which states that
an Artin algebra with infinite dominant dimension should be
self-injective. It is well-known that the finitistic dimension
conjecture implies the Nakayama conjecture.

Usually, a derived equivalence does not preserve the usual
homological dimensions of an algebra. However, under the condition
$\add(_AQ)=\add(\nu_AQ)$, we have the following inequalities about
these homological dimensions.

\begin{Koro}
Let $F:\Db{A}\longrightarrow\Db{B}$ be a derived equivalence between
Artin algebras $A$ and $B$. If $\add(_AQ)=\add(\nu_AQ)$, then

\smallskip
$(1)$  $\findim(A)\leq \findim(B)$,

$(2)$ $\gldim (A)\leq \gldim(B)$,

\smallskip
$(3)$  $\domdim(A)\geq \domdim(B)$. \label{ineqdimension}
\end{Koro}

{\it Proof.} $(1)$ For each $A$-module $X$, we have $\pd(_AX) =
\pd(_A\bar{G}\bar{F}(X))$ by Proposition \ref{fulldense}. According
to Proposition \ref{barGProp}(3), we have
$\pd(_A\bar{G}\bar{F}(X))\leq \pd(_B\bar{F}(X))$. By Proposition
\ref{barFProp}(3), we have another inequality $\pd(_B\bar{F}(X))\leq
\pd(_AX)$. Thus $\pd(_AX)=\pd(_B\bar{F}(X))$. This implies that
gl.dim$(A)\le $ gl.dim$(B)$. Moreover, if $\pd(_AX)<\infty$, we have
$\pd(_AX)=\pd(_B\bar{F}(X))\leq\findim(B)$. Hence $\findim(A)\leq
\findim(B)$. This proves $(1) and (2)$.

\smallskip
 $(3)$  Suppose $\domdim(B)=m$. Let $0\lra {}_BB\lra I_0\lra I_1\lra
\cdots $ be a minimal injective resolution of $_BB$. Then, by
definition, the injective $B$-modules $I_0,\cdots, I_{m-1}$ all are
projective. By Proposition \ref{barGProp}(1), we get an exact
sequence
$$0\lra {}_AQ^0\lra P_0\oplus \bar{G}(I_0)\lra P_1\oplus \bar{G}(I_1)\longrightarrow \cdots $$
with $P_i\in\add(_AQ)$. From this sequence we get another exact
sequence
$$0\lra {}_AQ^0\oplus {}_AQ\lra P_0\oplus \bar{G}(I_0)\oplus {}_AQ\lra P_1\oplus \bar{G}(I_1)\lra \cdots.$$
Since $\add(_AQ)=\add(\nu_AQ)$ and $\bar{G}(I_i)$ is injective for
all $i$ by Proposition \ref{barGProp}, the $A$-modules $P_i\oplus
\bar{G}(I_i)$ is injective for all $i$. Thus, the above exact
sequence actually gives an injective resolution of $_AQ^0\oplus
{}_AQ$. Set $I'_0:= P_0\oplus \bar{G}(I_0)\oplus {}_AQ$ and
$I'_i:=P_i\oplus \bar{G}(I_i)$ for $i>0$. Since $I_i$ is projective
for all $i\leq m-1$, we see that $\bar{G}(I_i)$ is projective for
all $0\le i\leq m-1$, and consequently $I'_i$ is projective for all
$0\le i\leq m-1$. Hence $\domdim(_AQ\oplus {}_AQ^0)\geq m$.
Moreover, since $\cpx{Q}$ is a tilting complex, we have
${}_AA\in\add(\bigoplus_{i}Q^i)$, that is, $_AA\in\add(_AQ\oplus
{}_AQ^0)$. Hence $\domdim(A)\geq m=\domdim(B)$. $\dickebox$

\section{Stable equivalences of Morita type induced by derived
equivalences\label{sectStM}}

In this section, we shall prove that an almost $\nu$-stable derived
functor $F$ between two finite-dimensional algebras actually induces
a stable equivalence of Morita type. Our result in this section
generalizes a well-known result of Rickard \cite[Corollary
5.5]{RickDFun}, which states that, for finite-dimensional
self-injective algebras, a derived equivalence induces a stable
equivalence of Morita type.

Throughout this section, we keep the notations introduced in Section
\ref{theFunctor} and consider exclusively finite-dimensional
algebras over a field $k$.

Let $\Lambda$ be an algebra. By $\Cz{\Lambda}$ (respectively,
$\Cf{\Lambda}$) we denote the full subcategory of $\C{\Lambda}$
consisting of all complexes bounded below (respectively, bounded
above). Analogously, one has the corresponding homotopy categories
$\Kz{\Lambda}$ and $\Kf{\Lambda}$ as well as the corresponding
derived categories $\Dz{\Lambda}$ and $\Df{\Lambda}$. Recall that
the category $\Df{\Lambda}$ is equivalent to the category
$\Kf{\pmodcat{\Lambda}}$, and the category $\Dz{\Lambda}$ is
equivalent to the category $\Kz{\imodcat{\Lambda}}$ (see
\cite[Theorem 10.4.8, p.388]{Weibel}, for example). Thus, for each
complex $\cpx{U}$ in $\Df{\Lambda}$ (respectively, $\Dz{\Lambda}$),
we can find a complex $\cpx{P_{U}}\in\Kf{\pmodcat{\Lambda}}$
(respectively, $\cpx{I_U}\in\Kz{\imodcat{\Lambda}}$) that is
isomorphic to $\cpx{U}$ in $\D{\Lambda}$.

Now, let $\cpx{X}$ be a complex in $\Df{\Lambda\opp}$ and $\cpx{Y}$
a complex in $\Df{\Lambda}$. By $\cpx{X}\otimesP_{\Lambda}\cpx{Y}$
we mean the total complex of the double complex with $(i,j)$-term
$X^i\otimes_{\Lambda}Y^j$, and by $\cpx{X}\otimesL_{\Lambda}\cpx{Y}$
the complex $\cpx{X}\otimesP_{\Lambda}\cpx{P_Y}$. Up to isomorphisms
of complexes, $\cpx{X}\otimesL_{\Lambda}\cpx{Y}$ does not depend on
the choice of $\cpx{P_Y}$. It is known that
$\cpx{X}\otimesL_{\Lambda}-$ is a functor from $\Df{\Lambda}$ to
$\Df{k}$, and is called the \emph{left derived functor} of
$\cpx{X}\otimesP_{\Lambda}-: \Kf{\Lambda}\longrightarrow\Kf{k}$.
Note that if we choose $\cpx{P}_X\in \Df{\Lambda\opp}$ such that
$\cpx{P}_X\simeq \cpx{X}$ in $\D{\Lambda\opp}$ then there is a
natural isomorphism between $\cpx{X}\otimesP_{\Lambda}\cpx{P_Y}$ and
$\cpx{P_X}\otimesP_{\Lambda}\cpx{Y}$ in $\D{k}$ (see \cite[Exercise
10.6.1, p.395]{Weibel}). Thus $\cpx{X}\otimesL_{\Lambda}\cpx{Y}$ can
be calculated by $\cpx{P_X}\otimesP_{\Lambda}\cpx{Y}$.

Let $\cpx{X_1}$ and $\cpx{X_2}$ be two complexes in $\Dz{\Lambda}$.
By $\HomP_{\Lambda}(\cpx{X_1},\cpx{X_2})$ we denote the total
complex of the double complex with $(i,j)$-term
$\Hom_{\Lambda}(X^{-i}_1, X^{j}_2)$. Choose
$\cpx{I_{X_2}}\in\Kz{\imodcat{\Lambda}}$ with $\cpx{I_{X_2}}\simeq
\cpx{X}_2$ in $\D{\Lambda}$. We define
$\rHom_{\Lambda}(\cpx{X_1},\cpx{X_2})= \HomP_{\Lambda}(\cpx{X_1},
\cpx{I_{X_2}}).$ It is known that $\rHom_{\Lambda}(\cpx{X_1}, -):
\Dz{\Lambda}\lra \Dz{k}$ is a functor. This functor is called the
\emph{right derived functor} of $\HomP_{\Lambda}(\cpx{X_1}, -):
\Kz{\Lambda}\lra\Kz{k}$. Note that if we choose
$\cpx{P}_{X_1}\in\Kf{\pmodcat{\Lambda}}$ with $\cpx{P}_{X_1}\simeq
\cpx{X}_1$ in $\D{\Lambda}$ then the complexes
$\HomP_{\Lambda}(\cpx{P_{X_1}}, \cpx{X_2})$ and
$\HomP_{\Lambda}(\cpx{X_1}, \cpx{I_{X_2}})$ are naturally isomorphic
in $\D{k}$ (see \cite[Exercise 10.7.1, p.400]{Weibel}). Thus
$\rHom_{\Lambda}(\cpx{X_1},\cpx{X_2})$ can be calculated by
$\HomP_{\Lambda}(\cpx{P_{X_1}}, \cpx{X_2})$.

Suppose $\Lambda_1$ and $\Lambda_2$ are two algebras. Let
$\cpx{T_i}$ be a tilting complex over $\Lambda_i$ with
$\Gamma_i=\End_{\Db{\Lambda_i}}(\cpx{T}_i)$ for $i=1,2.$ By a result
\cite[Theorem 3.1]{RickDFun} of Rickard,
$\cpx{T_1}\otimesP_k\cpx{T_2}$ is a tilting complex over
$\Lambda_1\otimes_k\Lambda_2$, and the endomorphism algebra of
$\cpx{T_1}\otimesP_k\cpx{T_2}$ is canonically isomorphic to
$\Gamma_1\otimes_k\Gamma_2$. Thus the tensor algebras
$\Lambda_1\otimes_k\Lambda_2$ and $\Gamma_1\otimes_k\Gamma_2$ are
derived-equivalent.

Recall that $\cpx{Q}$ is a tilting complex over $A$ with the
endomorphism algebra $B$. By \cite[Proposition 9.1]{RickMoritaTh},
$\HomP_A(\cpx{Q}, A)$ is a tilting complex over $A\opp$ with the
endomorphism algebra $B\opp$. Also, $\cpx{\bar{Q}}$ is a tilting
complex over $B$ with the endomorphism $A$, and therefore
$\HomP_B(\cpx{\bar{Q}}, B)$ is a tilting complex over $B\opp$ with
the endomorphism algebra $A$. Thus, by taking tensor products, we
get four derived-equivalent algebras $A\otimes_kA\opp$, $A\otimes_k
B\opp$, $B\otimes_k B\opp$, and $B\otimes_k A\opp$. The following
table, taken from \cite{RickDFun}, describes the corresponding
objects in various equivalent derived categories.

\renewcommand{\arraystretch}{1.4}
\begin{center}
\begin{tabular}{|c|c|c|c|}
\hline $\Db{A\otimes_kA\opp} $  & $\Db{A\otimes_k B\opp}$  &
$\Db{B\otimes_k B\opp}$
& $\Db{B\otimes_kA\opp}$\\
\hline $\cpx{Q}\otimesP_k A_A$ &
$\cpx{Q}\otimesP_k\HomP_B(\cpx{\bar{Q}}, B)$ &
$_BB\otimesP_k\HomP_B(\cpx{\bar{Q}},
B)$ & $_BB\otimes_kA_A$\\
$\cpx{Q}\otimesP_k\HomP_A(\cpx{Q},A)$ & $\cpx{Q}\otimesP_kB_B$
& $_BB\otimes_kB_B$ & $_BB\otimesP_k\HomP_A(\cpx{Q},A)$\\
$_AA\otimesP_k\HomP_A(\cpx{Q},A)$ & $_AA\otimes_kB_B$
& $\cpx{\bar{Q}}\otimesP_kB_B$ & $\cpx{\bar{Q}}\otimesP_k\HomP_A(\cpx{Q},A)$\\
$_AA\otimes_kA_A$ & $_AA\otimesP_k\HomP_B(\cpx{\bar{Q}}, B)$
& $\cpx{\bar{Q}}\otimesP_k\Hom_B(\cpx{\bar{Q}}, B)$ & $\cpx{\bar{Q}}\otimesP_kA_A$\\
$_AA_A$ & $\cpx{\Theta}$ & $_BB_B$ & $\cpx{\Delta}$\\
\hline
\end{tabular}

\medskip  Table 1
\end{center}

By Table 1, one can easily find the corresponding objects in the
above four equivalent derived categories. For instance, we consider
the derived equivalence $\widehat{F}:
\Db{A\otimes_kA\opp}\lra\Db{B\otimes_kA\opp}$ induced by the tilting
complex $\cpx{Q}\otimesP_kA_A$. Table 1 shows that $\widehat{F}$
sends $\cpx{Q}\otimesP_k A_A$ to $_BB\otimes_kA_A$,
$\cpx{Q}\otimesP_k\HomP_A(\cpx{Q},A)$ to
$_BB\otimesP_k\HomP_A(\cpx{Q},A)$, $_AA\otimesP_k\HomP_A(\cpx{Q},A)$
to $\cpx{\bar{Q}}\otimesP_kA_A$, and $_AA_A$ to $\cpx{\Delta}$. The
following lemma collects some properties of these complexes
\cite{RickDFun}.

\begin{Lem}
Let $\cpx{\Delta}$ and $\cpx{\Theta}$ be the complexes defined in
Table 1. We have the following.

$(1)$ $\cpx{\Delta}\otimesL_A\cpx{\Theta}\simeq {}_BB_B$ in
$\Db{B\otimes_kB\opp}$.

$(2)$ $\cpx{\Theta}\otimesL_B\cpx{\Delta}\simeq {}_AA_A$ in
$\Db{A\otimes_kA\opp}$.

$(3)$ The functor $\cpx{\Delta}\otimesL_A-: \Db{A}\lra \Db{B}$ is a
derived equivalence and $\cpx{\Delta}\otimesL_A\cpx{X}\simeq
F(\cpx{X})$ for all $\cpx{X}\in\Db{A}$.

$(4)$ The functor $\cpx{\Theta}\otimesL_B-:\Db{B}\lra\Db{A}$ is a
derived equivalence with $\cpx{\Theta}\otimesL_B\cpx{U}\simeq
G(\cpx{U})$ for all $\cpx{U}\in\Db{B}$.

$(5)$ $\cpx{\Theta}\simeq \rHom_B(\cpx{\Delta}, B)$ in
$\Db{A\otimes_kB\opp}$.

$(6)$ $\cpx{\Delta}$ is isomorphic to $\cpx{\bar{Q}}$ when
considered as an object in $\Db{B}$ and to $\Hom_A(\cpx{Q}, A)$ when
considered as an object in $\Db{A\opp}$.

$(7)$ $\cpx{\Theta}$ is isomorphic to $\cpx{{Q}}$ when considered as
an object in $\Db{A}$ and to $\Hom_B(\cpx{\bar{Q}}, B)$ when
considered as an object in $\Db{B\opp}$. \label{twotiltingproperty}
\end{Lem}

{\it Proof.} The statements (1)--(5) follow from \cite[Theorem 3.3,
Proposition 4.1]{RickDFun} and the remarks after \cite[Definition
4.2]{RickDFun}. The statements (6) and (7) are taken from
\cite[Proposition 3.1]{RickDFun}. $\dickebox$

\medskip
Note that it is an open question in \cite{RickDFun} whether the two
functors $F$ and $\cpx{\Delta}\otimesL_A-$ are naturally isomorphic,
although they have isomorphic images on each object by Lemma
\ref{twotiltingproperty}(3).

\medskip
Recall that a complex $\cpx{T}$ in $\Db{A\otimes_kB\opp}$ is called
a {\em two-sided tilting complex} over $A\otimes_kB\opp$ if there is
a complex $\cpx{\bar{T}}$ in $\Db{B\otimes_kA\opp}$ such that
$\cpx{T}\otimesL_{B}\cpx{\bar{T}}\simeq {}_AA_{A}$ in
$\Db{A\otimes_k A\opp}$ and $\cpx{\bar{T}}\otimesL_{A}\cpx{T}\simeq
{}_{B}B_{B}$ in $\Db{B\otimes_k B\opp}$. In this case, the complex
$\cpx{\bar{T}}$ is called an \emph{inverse} of $\cpx{T}$. From Lemma
\ref{twotiltingproperty} we see that $\cpx{\Delta}$ and
$\cpx{\Theta}$ defined in Table 1 are mutually inverse two-sided
tilting complexes over $A\otimes_kB\opp$ and $B\otimes_kA\opp$,
respectively.

\medskip
The following lemma, which is crucial in our later proofs, describes
some properties of the terms of the two-sided tilting complex
$\cpx{\Delta}$ in Table 1.

\begin{Lem}\label{deltasterms}
The two-sided tilting complex $\cpx{\Delta}$ is isomorphic in
$\Db{B\otimes_k A\opp}$ to a radical complex
$$0 \lra \Delta^0\lra \Delta^1\lra \dots \lra \Delta^n\lra 0 $$
with $\Delta^i\in\add( _B\bar{Q}\otimes_k Q^*_A)$ for all $i>0$.
\end{Lem}

{\it Proof.} Thanks to Table 1, there is a derived equivalence
$\widehat{F}: \Db{A\otimes_kA\opp}\longrightarrow
\Db{B\otimes_kA\opp}$. Moreover, the complexes
$\cpx{Q}\otimesP_kA_A$ and $\cpx{\bar{Q}}\otimesP_kA_A$ are the
associated tilting complexes to $\widehat{F}$ and its quasi-inverse,
respectively. Note that the two complexes are radical and have the
shape as assumed in Section \ref{theFunctor}. By Table 1, the
two-sided tilting complex $\cpx{\Delta}$ over $B\otimes_kA\opp$ is
isomorphic in $\Db{B\otimes_kA\opp}$ to $\widehat{F}(_AA_A)$, and
therefore, by Lemma \ref{formofFX}, the complex $\cpx{\Delta}$ is
isomorphic in $\Db{B\otimes_kA\opp}$ to a radical complex $\cpx{R}$:
$$0 \lra R^0\lra R^1\lra\cdots \lra R^n\lra 0$$ with
$R^i\in\add(_B\bar{Q}\otimes_kA_A)$ for all $i>0$.

Similarly, by Table 1, there is a derived equivalence
$\widetilde{F}:\Db{B\otimes_kB\opp}\longrightarrow\Db{B\otimes_kA\opp}$
induced by the tilting complex $_BB\otimesP_k\HomP_B(\cpx{\bar{Q}},
B)$ over $B\otimes_kB\opp$. From Table 1, we know that the complex
$_BB\otimesP_k\HomP_A(\cpx{Q},A)$ is a tilting complex associated to
the quasi-inverse of $\widetilde{F}$. Moreover, it follows from
Table 1 that $\widetilde{F}(_BB_B)\simeq\cpx{\Delta}$ in
$\Db{B\otimes_kA\opp}$. By Lemma \ref{formofFX},
$\widetilde{F}(_BB_B)$ is isomorphic in $\Db{B\otimes_kA\opp}$ to a
radical complex $\cpx{S}$:
$$0 \lra S^0\lra S^1\lra\cdots \lra S^n\lra 0,$$
with $S^i\in\add(_BB\otimes_kQ^*_A)$ for all $i>0$. Thus both
$\cpx{R}$ and $\cpx{S}$ are isomorphic to $\cpx{\Delta}$ in
$\Db{B\otimes_kA\opp}$. By Lemma \ref{kdiso}, the complexes
$\cpx{R}$ and $\cpx{S}$ are isomorphic in the homotopy category
$\Kb{B\otimes_kA\opp}$. Since $\cpx{R}$ and $\cpx{S}$ are radical
complexes, they are isomorphic in $\Cb{B\otimes_kA\opp}$. In
particular, $R^i\simeq S^i$ as $B$-$A$-bimodules for all $i$. Thus,
for each $i>0$, the bimodule $R^i$ lies in both
$\add(_BB\otimes_kQ^*_A)$ and $\add(_B\bar{Q}\otimes_kA_A)$. As a
result, we have $R^i\in\add(_B\bar{Q}\otimes_kQ^*_A)$ for all $i>0$.
$\dickebox$

\medskip
Using Lemma \ref{deltasterms}, we now prove the following main
result in this section.

\begin{Theo}
Let $A$ and $B$ be two finite-dimensional algebras over a field $k$,
and let $F:\Db{A}\lra\Db{B}$ be a derived equivalence. If $F$ is
almost $\nu$-stable, then there is a stable equivalence $\phi:
\stmodcat{A}\lra\stmodcat{B}$ of Morita type such that
$\phi(X)\simeq\bar{F}(X)$ for any $A$-module $X$, where $\bar{F}$ is
defined in Proposition \ref{stablefunctor}.
\label{derivedstbaleofMoritatype}
\end{Theo}

{\it Proof.} First, we show that $A$ may be assumed to be
indecomposable. In fact, if $A=A_1\times A_2$ is a product of two
algebras $A_i$, then the complex $\cpx{Q}$ associated to $F$ has a
decomposition $\cpx{Q}=\cpx{Q_1}\oplus\cpx{Q_2}$ such that
$\cpx{Q_1}\in\Kb{\pmodcat{A_1}}$ and
$\cpx{Q_2}\in\Kb{\pmodcat{A_2}}$. Correspondingly, the algebra $B$,
which is isomorphic to the endomorphism algebra of $\cpx{Q}$, is a
product of two algebras, say $B=B_1\times B_2$, such that $B_i\simeq
\End_{\Db{A_i}}(\cpx{Q}_i)$ for $i=1,2$. Thus the derived
equivalence $F:\Db{A}\lra \Db{B}$ induces two derived equivalences
$F_i: \Db{A_i}\lra\Db{B_i}$ for $i=1,2$. Moreover, for each $i$, the
complex $\cpx{Q_i}$ is a tilting complex associated to $F_i$, and
the tilting complex associated to the quasi-inverse of $F_i$ is
isomorphic to $F_i(A_i)\simeq F(A_i)$ which is a direct summand of
$\cpx{\bar{Q}}$. Thus, if $F$ is almost $\nu$-stable, then $F_i$ is
almost $\nu$-stable for $i=1, 2$. Furthermore, if $A_i$ and $B_i$
are stably equivalent of Morita type for $i=1,2$, then $A_1\times
A_2$ and $B_1\times B_2$ are stably equivalent of Morita type. Thus,
we may assume that $A$ is indecomposable.

Since derived equivalence preserves the semisimplicity of algebras,
we know that $A$ is semi-simple if and only if $B$ is semisimple.
Hence we can further assume that $A$ is non-semisimple. Now, let $A$
be non-semisimple and indecomposable. Then $B$ is also
non-semisimple and indecomposable.

Let $\cpx{\Delta}$ be the complex in Table 1. By Lemma
\ref{deltasterms}, the complex $\cpx{\Delta}$ is isomorphic in
$\Db{B\otimes_kA\opp}$ to a radical complex
$$
0\lra \Delta^0\lra \Delta^1\lra \cdots\lra \Delta^n\lra 0,
$$ with $\Delta^i$ in $\add(_B\bar{Q}\otimes_kQ^*_A)$ for all $i>0$.
For simplicity, we assume that $\cpx{\Delta}$ is the above complex.
By Lemma \ref{twotiltingproperty} (6), the complex $\cpx{\bar{Q}}$
is isomorphic to $\cpx{\Delta}$ in $\Db{B}$, and therefore in
$\Kb{B}$ by Lemma \ref{kdiso}(1). Thus there is a chain map $\alpha:
\cpx{\bar{Q}}\rightarrow\cpx{\Delta}$ such that the mapping cone
$\con(\alpha)$ is acyclic. Since the terms of $\con(\alpha)$ in
positive degrees are all projective, the acyclic complex
$\con(\alpha)$ splits. This means that $_B\Delta^0$ is a projective
$B$-module. Thus all terms of $\cpx{\Delta}$ are projective as left
$B$-modules. Similarly, by the fact that $\cpx{\Delta}$ is
isomorphic to $\Hom_A(\cpx{Q}, A)$ in $\Db{A\opp}$, we infer that
$\Delta^0_A$ is a projective right $A$-module and that all the terms
of $\cpx{\Delta}$ are projective as right $A$-modules. Consequently,
the complex $\rHom_B(\cpx{\Delta}, B)$ is isomorphic in
$\Db{A\otimes_kB\opp}$ to the complex $\HomP_B(\cpx{\Delta}, B)$:
$$\xymatrix{
0\ar[r] & \Hom_B(\Delta^n, B)\ar[r]& \cdots\ar[r] & \Hom_B(\Delta^1,
B)\ar[r]  & \Hom_B(\Delta^0, B)\ar[r] & 0.
 }$$
Since $\Delta^i\in\add (_B\bar{Q}\otimes_k Q^*_A)$ for all positive
$i$, we find that the $A$-$B$-bimodule $\Hom_B(\Delta^i, B)$ is in
$\add (\Hom_B(_B\bar{Q}\otimes_k Q^*_A, B))$ for each positive
integer $i$. Recall that there is a natural isomorphism
$$\Hom_B(_B\bar{Q}\otimes_k Q^*_A, B)\simeq \nu_A
Q\otimes_k\bar{Q}^*.$$ Since $\add(_AQ)=\add(\nu_AQ)$, the term
$\Hom_B(\Delta^i, B)$ is in $\add (Q\otimes\bar{Q}^*)$ for all
$i>0$. As to $\Hom_B(_B\Delta^0_A, B)$, we can use Lemma
\ref{twotiltingproperty} (7) and similarly prove that
$\Hom_B(_B\Delta^0_A, B)$ is projective as a one-sided module on
both sides.

Next, we show that $\Delta^0$ and $\Hom_B(\Delta^0, B)$ define a
stable equivalence of Morita type between $A$ and $B$. Indeed, by
Lemma \ref{twotiltingproperty} (5), the two-sided tilting complex
$\cpx{\Theta}$ defined in Table 1 is isomorphic in
$\Db{A\otimes_kB\opp}$ to $\rHom_B(\cpx{\Delta}, B)$, and the latter
is isomorphic in $\Db{A\otimes_kB\opp}$ to $\HomP_B(\cpx{\Delta},
B)$. For simplicity, we assume that $\cpx{\Theta}$ equals
$\HomP_B(\cpx{\Delta}, B)$. Since all the terms of $\cpx{\Delta}$
are projective as right $A$-modules, the complex
$\cpx{\Delta}\otimesL_A\cpx{\Theta}$ is isomorphic in
$\Db{B\otimes_kB\opp}$ to the complex
$\cpx{\Delta}\otimesP_A\cpx{\Theta}$. The $m$-th term of
$\cpx{\Delta}\otimesP_A\cpx{\Theta}$ is
$$\bigoplus_{i+j=m}(\Delta^i\otimes_A\Theta^j)=\bigoplus_{i+j=m}(\Delta^i\otimes_A\Hom_B(\Delta^{-j},
B)).$$ Let $_AE$ be the maximal $\nu$-stable $A$-module, and let
$_B\bar{E}$ be the maximal $\nu$-stable $B$-module (see Section
\ref{theFunctor}). Since $\add(_AQ)=\add(\nu_AQ)$ and
$\add(_B\bar{Q})=\add(\nu_B\bar{Q})$, we have $_AQ\in\add(_AE)$ and
$_B\bar{Q}\in\add(_B\bar{E})$. By Lemma \ref{twotiltingproperty} (3)
and  Lemma \ref{lem4.2}, the complex $\cpx{\Delta}\otimesP_AE$ is
isomorphic in $\Db{B}$ to a complex in $\Kb{\add(_B\bar{E})}$. Since
$\Delta^i\in\add(_B\bar{Q}\otimes_kQ_A)$ for $i>0$, the $i$-th term
$\Delta^i\otimes_AE$ of $\cpx{\Delta}\otimesP_AE$ is in
$\add(_B\bar{Q})\subseteq\add(_B\bar{E})$ for all $i>0$. It follows
that $\Delta^0\otimes_AE$ is in $\add(_B\bar{E})$, and therefore
$$E^*\otimes_A\Theta^0=E^*\otimes_A\Hom_B(\Delta^0, B)\simeq\Hom_B(\Delta^0\otimes_AE,
B)\in\add(\bar{E}^*).$$ Note that, for each $i>0$, we have
$\Delta^i\in\add(\bar{Q}\otimes_kQ^*)\subseteq\add(\bar{E}\otimes_kE^*)$
and $\Theta^{-i}=\Hom_B(\Delta^i,
B)\in\add(Q\otimes_k\bar{Q}^*)\subseteq\add(E\otimes_k\bar{E}^*)$.
Thus, it is not hard to see that $\Delta^i\otimes_A\Theta^j$ is in
$\add(\bar{E}\otimes_k\bar{E}^*)$ for all $i$ and $j$ with $i+j\neq
0$. This means that all terms in non-zero degrees of
$\cpx{\Delta}\otimesP_A\cpx{\Theta}$ are in
$\add(\bar{E}\otimes_k\bar{E}^*)$. For the term in degree $0$ of
$\cpx{\Delta}\otimesP_A\cpx{\Theta}$, except
$\Delta^0\otimes_A\Theta^0$, all of its other direct summands are in
$\add(\bar{E}\otimes_kE^*\otimes_AE\otimes_k\bar{E}^*)$ which is
contained in $\add(\bar{E}\otimes_k\bar{E}^*)$. Note that all the
bimodules in $\add(_B\bar{E}\otimes_k\bar{E}^*_B)$ are
projective-injective.

Now, we have $\cpx{\Delta}\otimesP_A\cpx{\Theta}\simeq {}_BB_B$ in
$\Db{B\otimes_kB\opp}$ by Lemma \ref{twotiltingproperty} (1). Thus
the complex $\cpx{\Delta}\otimesP_A\cpx{\Theta}$ has zero homology
and projective-injective terms in all non-zero degrees, and
therefore it splits and is isomorphic to $_BB_B$ in the homotopy
category $\Kb{B\otimes_kB\opp}$. Since $B$ is indecomposable and
non-semisimple, the bimodule ${}_BB_B$ is indecomposable and
non-projective, and therefore it is a direct summand of
$\Delta^0\otimes_A\Theta^0$. It follows that
$\Delta^0\otimes_A\Theta^0\simeq {}_BB_B\oplus U$ for a
projective-injective $B$-$B$-bimodule $U$. Similarly, we have
$\Theta^0\otimes_B\Delta^0\simeq {}_AA_A\oplus V$ for a
projective-injective $A$-$A$-bimodule $V$. Hence $A$ and $B$ are
stably equivalent of Morita type.

Let $\phi: \stmodcat{A}\lra\stmodcat{B}$ be the stable equivalence
induced by $\Delta^0\otimes_A-$. It follows from Lemma
\ref{twotiltingproperty} $(3)$ and Lemma \ref{formofFX} that
$\phi(X)\simeq \bar{F}(X)$ in $\stmodcat{B}$ for all $A$-modules
$X$. $\dickebox$

\medskip
Let $A$ be an algebra. An $A$-module $M$ is called a
\emph{generator-cogenerator} for $A$-mod if $A\oplus
D(A)\in\add(M)$. The {\em representation dimension} of $A$, denoted
by $\repdim (A)$, is defined to be
$$\repdim(A):=\inf\{\gldim(\End_A(M))\mid M \mbox{ is a generator-cogenerator for}\;A\mbox{-mod}\}.$$
This notion was introduced by Auslander in \cite{A} to measure
homologically how far an algebra is from being
representation-finite, and has been studied by many authors in
recent years (see \cite{Rouq} and the references therein).

The following result is a consequence of Theorem
\ref{derivedstbaleofMoritatype} since stable equivalences of Morita
type preserve representation dimensions \cite{Xi1}.

\begin{Koro} If $F$ is an almost $\nu$-stable derived
functor between $A$ and $B$, then $A$ and $B$ have the same
representation and dominant dimensions.
\end{Koro}

\medskip
As another consequence of Theorem \ref{derivedstbaleofMoritatype},
we re-obtain the following result of Rickard \cite{RickDstable}
since every derived equivalence between self-injective algebras is
almost $\nu$-stable by Proposition \ref{prop4.1}.

\begin{Koro} Let $A$ and $B$ be finite-dimensional self-injective algebras. If $A$ and
$B$ are derived-equivalent, then they are stably equivalent of
Morita type.
\end{Koro}

{\it Remark:} (1) Let $A$ be a finite dimensional self-injective
algebra and $X$ be an $A$-module. By the remark at the end of
Section \ref{theFunctor}, there is a derived equivalence between the
algebras $\End_A(A\oplus X)$ and $\End_A(A\oplus\Omega(X))$
satisfying the almost $\nu$-stable condition. Thus, we have an
alternative proof of the result \cite[Corollary 1.2]{LiuXi3} of Liu
and Xi by applying Theorem \ref{derivedstbaleofMoritatype}.

(2) Theorem \ref{derivedstbaleofMoritatype} may be false if only one
of the two equalities of the almost $\nu$-stable condition is
satisfied. For a counterexample, we refer the reader to Example 2 in
Section \ref{example}.


\section{Inductive constructions of almost $\nu$-stable derived equivalences}\label{sectConst}

In this section, we shall give several inductive constructions of
almost $\nu$-stable derived equivalences. As a consequence, one can
produce a lot of (usually not self-injective) finite-dimensional
algebras that are both derived-equivalent and stably equivalent of
Morita type.

In this section, we keep the notations introduced in Section
\ref{theFunctor}. Our first inductive construction is the following
proposition.

\begin{Prop}\label{EndDerStable}
Suppose that $F$ is an almost $\nu$-stable derived equivalence
between finite-dimensional algebras $A$ and $B$ over a field $k$.
Let $\bar{F}$ be the stable functor of $F$ defined in Proposition
\ref{stablefunctor}, and let $X$ be an $A$-module. Then there is an
almost $\nu$-stable derived equivalence between the endomorphism
algebras $\End_A(A\oplus X)$ and $\End_B(B\oplus\bar{F}(X))$.
\end{Prop}

{\it Proof.} We keep the notations in the proof of Theorem
\ref{derivedstbaleofMoritatype}. By the last part of the proof of
Theorem \ref{derivedstbaleofMoritatype}, the two-sided tilting
complexes $\cpx{\Delta}$ and $\cpx{\Theta}$ have the properties:
$\cpx{\Delta}\otimesP_A\cpx{\Theta}\simeq {}_BB_B$ in
$\Kb{B\otimes_k B\opp}$ and
$\cpx{\Theta}\otimesP_B\cpx{\Delta}\simeq {}_AA_A$ in
$\Kb{A\otimes_k A\opp}$. It follows that the functor
$(\cpx{\Delta}\otimesP_A\cpx{\Theta})\otimesP_B-$ is naturally
isomorphic to the identity functor $1_{\Kb{B}}$, and
$(\cpx{\Theta}\otimesP_B\cpx{\Delta})\otimesP_A-$ is naturally
isomorphic to the identity functor $1_{\Kb{A}}$. Thus
$\cpx{\Delta}\otimesP_A-$ and $\cpx{\Theta}\otimesP_B-$  induce
mutually inverse equivalences between $\Kb{A}$ and $\Kb{B}$. Now we
prove that the restrictions of these two functors to
$\Kb{\add(A\oplus X)}$ and to $\Kb{\add(B\oplus\bar{F}(X))}$ are
also mutually inverse equivalences for each $A$-module $X$.

In fact, the complex $\cpx{\Delta}\otimesP_AX$ is of the following
form
$$\xymatrix{
0\ar[r] & \Delta^0\otimes_AX\ar[r] & \Delta^1\otimes_A X\ar[r] &
\cdots\ar[r] & \Delta^n\otimes_AX\ar[r] & 0.
 }$$
Since $\Delta^i$ is a projective bimodule for all $i>0$, the term
$\Delta^i\otimes_AX$ is a projective $B$-module for all $i>0$.
Moreover, by Theorem \ref{derivedstbaleofMoritatype}, we have
$\Delta^0\otimes_AX\simeq\bar{F}(X)$ in $\stmodcat{B}$, and
therefore $\Delta^0\otimes_AX$ is a direct summand of
$\bar{F}(X)\oplus P$ for some projective $B$-module $P$. Hence the
complex $\cpx{\Delta}\otimesP_AX$ is in
$\Kb{\add(B\oplus\bar{F}(X))}$. Note that, for each projective
$A$-module $P_1$, the complex $\cpx{\Delta}\otimesP_AP_1$ is in
$\add(\cpx{\bar{Q}})$. Thus, for each complex $\cpx{X}$ in
$\Kb{\add(A\oplus X)}$, the complex $\cpx{\Delta}\otimesP_A\cpx{X}$
is in $\Kb{\add(B\oplus\bar{F}(X))}$. Similarly, the functor
$\cpx{\Theta}\otimesP_B-$ takes complexes in
$\Kb{\add(B\oplus\bar{F}(X))}$ to complexes in $\Kb{\add(A\oplus
X)}$. Thus $\cpx{\Delta}\otimesP_A-$ and $\cpx{\Theta}\otimesP_B-$
induce mutually inverse equivalences between the triangulated
categories $\Kb{\add(A\oplus X)}$ and
$\Kb{\add(B\oplus\bar{F}(X))}$.

Let $\Lambda=\End_A(A\oplus X)$ and
$\Gamma=\End_B(B\oplus\bar{F}(X))$. Then $\Kb{\pmodcat{\Lambda}}$
and $\Kb{\pmodcat{\Gamma}}$ are canonically equivalent to
$\Kb{\add(A\oplus X)}$ and $\Kb{\add(B\oplus\bar{F}(X))}$,
respectively. By \cite[Theorem 6.4]{RickMoritaTh}, there is a
derived equivalence $\widehat F$ between $\Lambda$ and $\Gamma$.
Moreover, the tilting complexes associated to $\widehat F$ and its
quasi-inverse are $\HomP_A(A\oplus
X,\cpx{Q}\oplus(\cpx{\Theta}\otimesP_B\bar{F}(X)))$ and
$\HomP_B(B\oplus\bar{F}(X),
\cpx{\bar{Q}}\oplus(\cpx{\Delta}\otimesP_AX))$, respectively. By the
proof of Theorem \ref{derivedstbaleofMoritatype}, the $i$-th term
$\Theta^i$ of $\cpx{\Theta}$ is in $\add(Q\otimes_k\bar{Q}^*)$ for
all $i<0$. Hence $\Theta^i\otimes_B\bar{F}(X)$ is in $\add(Q)$ for
all $i<0$, and all the terms in negative degrees of $\HomP_A(A\oplus
X,\cpx{Q}\oplus(\cpx{\Theta}\otimesP_B\bar{F}(X)))$  are in
$\add(\Hom_A(A\oplus X, Q))$. Similarly, all the terms in positive
degrees of $\HomP_B(B\oplus\bar{F}(X),
\cpx{\bar{Q}}\oplus(\cpx{\Delta}\otimesP_AX))$ are in
$\add(\Hom_B(B\oplus\bar{F}(X), \bar{Q})$. Note that we have the
following isomorphisms
$$\begin{array}{rl}
\nu_{\Lambda}(\Hom_A(A\oplus X, Q)) & =
D\Hom_{\Lambda}(\Hom_A(A\oplus X, Q),\Hom_A(A\oplus X, A\oplus X))\\
& \simeq D\Hom_A(Q, A\oplus X)\\
&\simeq D(\Hom_A(Q, A)\otimes_A(A\oplus X))\\
& \simeq \Hom_A(A\oplus X, D(\Hom_A(Q, A)))\\
& = \Hom_A(A\oplus X, \nu_AQ).
\end{array}$$
Since $\add(_AQ)=\add(\nu_AQ)$, we have $\add(\Hom_A(A\oplus X,
Q))=\add(\nu_{\Lambda}(\Hom_A(A\oplus X, Q)))$. Similarly, we have
$\add(\Hom_B(B\oplus \bar{F}(X),
\bar{Q}))=\add(\nu_{\Gamma}(\Hom_B(B\oplus\bar{F}(X), \bar{Q})))$.
This shows that the derived equivalence between $\Lambda$ and
$\Gamma$ induced by the tilting complex $\HomP_A(A\oplus
X,\cpx{Q}\oplus(\cpx{\Theta}\otimesP_B\bar{F}(X)))$ is almost
$\nu$-stable. $\dickebox$

\medskip
Our next construction uses tensor products.

\begin{Prop}
Let $k$ be a field. Suppose $F$ is an almost $\nu$-stable  derived
equivalence between finite-dimensional $k$-algebras $A$ and $B$.
Then, for each finite-dimensional self-injective $k$-algebra $C$,
there is an almost $\nu$-stable derived equivalence between the
tensor algebras $A\otimes_k C\opp$ and $B\otimes_kC\opp$.
\end{Prop}

{\it Proof.} By \cite[Theorem 2.1]{RickDFun}, $F$ induces a derived
equivalence $\widehat{F}$ between $A\otimes_k C\opp$ and
$B\otimes_kC\opp$. Suppose $\cpx{Q}$ and $\cpx{\bar{Q}}$ be the
associated tilting complexes of $F$ and its quasi-inverse,
respectively. From Table 1 we know that $\cpx{Q}\otimesP_k C_C$ and
$\cpx{\bar{Q}}\otimesP_k C_C$ are the associated tilting complexes
of $\widehat{F}$ and its quasi-inverse, respectively. Now we have
the following isomorphisms:
$$\begin{array}{rcl}
\nu_{A\otimes_kC\opp}(_AQ\otimes_kC_C)& = &D\Hom_{A\otimes_kC\opp}(_AQ\otimes_kC_C, {}_AA\otimes_k C_C)\\
& \simeq & D\big(\Hom_A(_AQ, {}_AA)\otimes_k\Hom_{C\opp}(C_C, C_C)\big)\\
& \simeq & D\Hom_A(_AQ, {}_AA)\otimes_kD\Hom_{C\opp}(C_C, C_C)\\
& \simeq & \nu_A Q\otimes_k\nu_{C^{\opp}}C_C \\
& \simeq & \nu_A Q\otimes_kC_C \quad (\mbox{because\;} C \mbox{\; is
self-injective}).
\end{array}$$
Since $F$ is almost $\nu$-stable, we have $\add(Q)=\add(\nu_AQ)$ and
$\add(\nu_{A\otimes_kC\opp}({}_AQ\otimes_kC_C))=\add({}_AQ\otimes_kC_C)$.
Similarly, we have
$\add(\nu_{B\otimes_kC\opp}(_B\bar{Q}\otimes_kC_C))=\add(_B\bar{Q}\otimes_kC_C)$.
Hence $\widehat{F}$ is almost $\nu$-stable and the proof is
completed. $\dickebox$

\medskip
Let $A$ be a finite-dimensional algebra over a field $k$, and let
$X$ be an $A$-module. The {\em one-point extension} of $A$ by $X$,
denoted by $A[X]$, is the triangular matrix algebra
$\left[{{k}\atop{X}}\,{{0}\atop{A}}\right]$. The natural projection
$A[X]\longrightarrow A$ shows that $A$ is a quotient algebra of
$A[X]$, and that $\modcat{A}$ can be viewed as a full subcategory of
$\modcat{A[X]}$. Let $\widetilde{X}$ denote the $A[X]$-module
$\left[{k\atop X}\right]$. Then, for each $A$-module $M$, we have
$\Hom_{A[X]}(M, \widetilde{X})\simeq\Hom_A(M, X)$.

\smallskip
Our third construction of an almost $\nu$-stable derived equivalence
is given by one-point extensions.

\begin{Prop}
Let $k$ be field. Suppose $F$ is an almost $\nu$-stable derived
equivalence between finite-dimensional $k$-algebras $A$ and $B$. If
$X$ is an $A$-module such that $F(X)$ is isomorphic in $\Db(B)$ to a
$B$-module $Y$, then there is an almost $\nu$-stable derived
equivalence between the one-point extensions $A[X]$ and $B[Y]$.
\label{1point}
\end{Prop}

{\it Proof.} Let $G$ be a quasi-inverse of $F$. Recall that
$\cpx{Q}$ and $\cpx{\bar{Q}}$ denote the radical tilting complexes
associated to $F$ and $G$, respectively. Then $\cpx{Q}$ can be
viewed as a complex in $\Kb{\pmodcat{A[X]}}$. By a result of Barot
and Lenzing \cite{OnepoitExtension}, the complex
$\cpx{Q}\oplus\widetilde{X}$ is a tilting complex over $A[X]$ such
that its endomorphism algebra is isomorphic to $B[Y]$, where
$\widetilde{X}$ is regarded as a complex concentrated only on degree
zero. Moreover, $\cpx{\bar{Q}}\oplus\widetilde{Y}$ is a tilting
complex associated to the quasi-inverse of the derived equivalence
induced by $\cpx{Q}\oplus\widetilde{X}$.  Recall that $Q$ is the
direct sum of all the terms of $\cpx{Q}$ in negative degrees. Then
$\add (\nu_A Q)=\add(_AQ)$ by assumption. Since the direct sum of
all terms in negative degrees of $\cpx{Q}\oplus\widetilde{X}$ equals
$Q$, we have to show that $\add(\nu_{A[X]}Q)=\add(_{A[X]}Q)$.

Since $F(X)$ is isomorphic in $\Db{B}$ to the $B$-module $Y$, we
have $\Hom(\cpx{Q}, X[i])=0$ for all $i\neq 0$. Then there is a
unique maximal submodule $L$ of $X$ with respect to the property
$\Hom_A(Q, L)=0$. This shows that $0=\Hom_{\Db{A}}(\cpx{Q},
X[i])\simeq\Hom_{\Db{A}}(\cpx{Q}, (X/L)[i])$ for all integers $i>0$.
If $X/L\neq 0$, then $\Hom_A(Q, \soc(X/L))$ $\neq 0$ by the
definition of $L$. This implies that $\Hom_{\Db{A}}(\cpx{Q},
(X/L)[i])\neq 0$ for some $i>0$, a contradiction. Thus $X/L=0$,
$\Hom_A(Q, X)=0$, and $\Hom_{A[X]}(Q, \widetilde{X})=0$.
Consequently,
$$\nu_{A[X]}Q=D\Hom_{A[X]}(Q, A[X])\simeq D\Hom_{A[X]}(Q, A\oplus\widetilde{X})
= D\Hom_{A[X]}(Q, A)\simeq \nu_A Q.$$ Hence
$\add(\nu_{A[X]}(Q))=\add(_{A[X]}Q)$. Similarly, we have
$\add(\nu_{B[Y]}\bar{Q})=\add(_{B[Y]}\bar{Q})$. This finishes the
proof. $\dickebox$

\section{Examples and questions}\label{example}

In the following, we shall illustrate our results with examples.

{\parindent=0pt\bf Example 1}: Let $A$ and $B$ be finite-dimensional
$k$ algebras given by quiver with relations in Fig. 1 and Fig. 2,
respectively.
\begin{center}
\begin{tabular}{ccc}
\multirow{2}{*}{$\xymatrix@R=2.4mm@C=3mm{
 \bullet\ar[rr]^{\alpha}^(0){1}^(1){2} && \bullet\ar[ldd]^{\beta}\\
      &&\\
      &\bullet\ar[luu]^(.15){3}^{\gamma} &
}$} & \hspace{2cm} & $\xymatrix{
  \bullet\ar@<2.5pt>[r]^{\alpha} &
  \bullet\ar@<2.5pt>[l]^(1){1}^(0){2}^{\beta}\ar@<2.5pt>[r]^{\gamma}
  &\bullet\ar@<2.5pt>[l]^(0){3}^{\delta}
}$\\
&& $\alpha\gamma=\delta\beta=0$\\
$\alpha\beta\gamma\alpha=\beta\gamma\alpha\beta=\gamma\alpha\beta\gamma=0$
&&
$\alpha\beta\alpha=\delta\gamma\delta=\beta\alpha-\gamma\delta=0.$\\
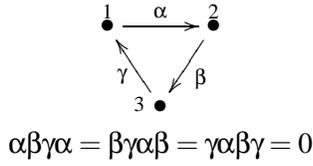
{\footnotesize Fig. 1 } && {\footnotesize Fig. 2 }\\
\end{tabular}
\end{center}
Let $P_A(i)$, $I_A(i)$ and $S_A(i)$ denote the indecomposable
projective, injective and simple $A$-modules corresponding to the
vertex $i$, respectively. We take a non-zero homomorphism $f:
P_A(2)\ra P_A(1)$. Then there is a tilting complex of $A$-modules
$$\cpx{Q}: \quad \xymatrix{0\ar[r] & P_A(2)\oplus P_A(2)\oplus P_A(3)
\ar[r]^(.70){[f,0,0]^T} & P_A(1)\ar[r] & 0.}$$ The endomorphism
algebra of $\cpx{Q}$ is isomorphic to $B$.  Let $F:
\Db{A}\rightarrow\Db{B}$ be a derived equivalence with $\cpx{Q}$ as
its associated tilting complex. Clearly, $F$ is almost $\nu$-stable
since $A$ and $B$ are symmetric algebras. By \cite[Proposition
7.3]{HuXi2}, there is a derived equivalence $F_1$ between
$\bar{A}=A/(\alpha\beta\gamma)$ and $\bar{B}=B/(\alpha\beta)$ with
inverse $G_1$ such that the associated tilting complexes over
$\bar{A}$ and over $\bar{B}$ are
$$\cpx{Q}_1: \quad \xymatrix{0\ar[r] & P_{\bar{A}}(2)\oplus P_{\bar{A}}(2)
\oplus P_{\bar{A}}(3)\ar[r] & P_A(1)/\soc P_A(1)\ar[r] & 0,}$$
$$\cpx{\bar{Q}}_1: \quad \xymatrix{0\ar[r] &P_B(1)/\soc P_B(1) \ar[r]
& P_{\bar{B}}(2)\oplus P_{\bar{B}}(2)\oplus P_{\bar{B}}(3)\ar[r] &
0,}$$ respectively. Clearly, the two complexes satisfy the
conditions: $\add(_{\bar{A}}Q_1)=\add(\nu_{\bar{A}}Q_1)$ and
$\add(\bar{Q}_1)=\add(\nu_{\bar{B}}\bar{Q}_1)$. Hence the algebras
$\bar{A}$ and $\bar{B}$ are both derived-equivalent and stably
equivalent of Morita type by Theorem
\ref{derivedstbaleofMoritatype}. We know that $F_1(S_{\bar{A}}(1))$
is isomorphic to the simple $\bar{B}$-module $S_{\bar{B}}(1)$. The
one-point extension $\bar{A}[S_{\bar{A}}(1)]$ is given by the quiver
Fig. 3.
$$\begin{array}{ccc}
 \xymatrix@R=2.4mm@C=3mm{
  &\bullet\ar[ld]_(-0.20){4}_{\eta}&\\
 \bullet\ar[rr]_{\alpha}^(0){1}^(1){2} && \bullet\ar[ldd]^{\beta}\\
      &&\\
      &\bullet\ar[luu]^(-0.20){3}^{\gamma} &
}&\quad\quad\quad\quad& \xymatrix{
  \bullet\ar[d]^(0){4}^{\eta}&&\\
  \bullet\ar@<2.5pt>[r]^{\alpha} &
  \bullet\ar@<2.5pt>[l]^(1){1}^(0){2}^{\beta}\ar@<2.5pt>[r]^{\gamma}
  &\bullet\ar@<2.5pt>[l]^(0){3}^{\delta}
}\\
\mbox{\rm\footnotesize Fig. 3} &\quad\quad\quad\quad&
\mbox{\rm\footnotesize Fig. 4}
\end{array}$$
with relations
$\alpha\beta\gamma=\beta\gamma\alpha\beta=\gamma\alpha\beta\gamma=\eta\alpha=0$,
and the one-point extension $\bar{B}[S_{\bar{B}}(1)]$ is given by
the quiver Fig. 4. with relations
$\eta\alpha=\alpha\beta=\delta\gamma\delta=\beta\alpha-\gamma\delta=\alpha\gamma=\delta\beta=0$.
By Proposition \ref{1point}, there is a derived equivalence between
$\bar{A}[S_{\bar{A}}(1)]$ and $\bar{B}[S_{\bar{B}}(1)]$, which
induces a stable equivalence of Morita type.

An calculation shows that ${F}_1(I_{\bar{A}}(1))$ is isomorphic to
the $\bar{B}$-module $I_{\bar{B}}(1)$. The algebras
$\End_{\bar{A}}(\bar{A}\oplus I_{\bar{A}}(1))$ and
$\End_{\bar{B}}(\bar{B}\oplus I_{\bar{B}}(1))$ are given by Fig. 5
and Fig. 6, respectively.
$$\begin{array}{ccc}
\xymatrix@!=5mm{
\bullet\ar[r]^(0){1}^(1){2}^{\alpha} & \bullet\ar[d]^{\beta}\\
\bullet\ar[u]^{\delta} & \bullet\ar[l]^{\gamma}^(0){3}^(1){4}
}&\quad\quad\quad & \xymatrix{
& \bullet\ar[rd]^{\eta} &\\
\bullet\ar@<2.5pt>[r]^{\alpha} &
\bullet\ar@<2.5pt>[l]^(1){1}^(0){2}^{\beta}\ar[u]^(1){4}^{\delta}
&\bullet\ar[l]^{\gamma}^(0){3}
}\\
\alpha\beta\gamma\delta=\beta\gamma\delta\alpha=\gamma\delta\alpha\beta\gamma=0
&
&\alpha\delta=\gamma\beta=\beta\alpha-\delta\eta\gamma=\gamma\delta\eta=\eta\gamma\delta=0\\
\mbox{\rm\footnotesize Fig. 5 }&& \mbox{\rm\footnotesize Fig. 6}
\end{array}$$
Thus $\End_{\bar{A}}(\bar{A}\oplus I_{\bar{A}}(1))$ and
$\End_{\bar{B}}(\bar{B}\oplus I_{\bar{B}}(1))$ are
derived-equivalent and stably equivalent of Morita type by
Proposition \ref{EndDerStable}.

\medskip
The following example, taken from \cite{RickDstable}, shows that
Theorem \ref{derivedstbaleofMoritatype} may fail if only one of the
conditions of an almost $\nu$-stable functor is satisfied.

{{\parindent=0pt\bf Example 2:} Let $A$ be the $17$-dimensional
algebra given by the quiver
$$
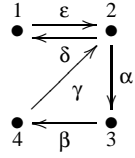
\xymatrix{
 \bullet\ar@<2.5pt>[r]^{\epsilon}^(0){1}^(1){2} & \bullet\ar@<2.5pt>[l]^{\delta}\ar[d]^{\alpha}\\
 \bullet\ar[ru]_{\gamma} & \bullet\ar[l]^{\beta}^(0){3}^(1){4}\\
 }$$
with relations $\gamma\alpha\beta=\gamma\delta=
\epsilon\alpha\beta=0, \delta\epsilon=\alpha\beta\gamma$. As before,
we denote by $P_A(i)$ the indecomposable projective $A$-module
corresponding to the vertex $i$. Let $\cpx{Q}$ be the direct sum of
the following two complexes
$$\begin{array}{ccccccc} 0 & \lra & P_A(1)                      & \lra & P_A(2)& \lra  & 0,\\
                   0 & \lra & 0 & \lra & P_A(2)\oplus P_A(3)\oplus P_A(4)   &  \lra & 0, \end{array}$$
where $P_A(1)$ is in degree $-1$. One can check that $\cpx{Q}$ is a
tilting complex over $A$. Let $B=\End_{\Db{A}}(\cpx{Q})$. Then $B$
is a 20-dimensional algebra given by the quiver
$$
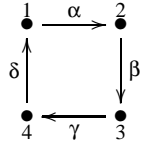
\xymatrix{
\bullet\ar[r]^{\alpha}^(0){1}^(1){2} & \bullet\ar[d]^{\beta}\\
\bullet\ar[u]^{\delta} & \bullet\ar[l]^{\gamma}^(0){3}^(1){4} }$$
with relations
$\alpha\beta\gamma\delta\alpha=0=\delta\alpha\beta\gamma$, where the
indecomposable projective $B$-modules at the vertices $1$, $2$, $3$,
and $4$ correspond respectively to the direct summands
$P_A(1)\rightarrow P_A(2)$, $P_A(2)$, $P_A(3)$ and $P_A(4)$ of the
complex $\cpx{Q}$. Let $F: \Db{A}\longrightarrow\Db{B}$ be the
derived equivalence induced by the tilting complex $\cpx{Q}$. Then
$F(P_A(i))=P_B(i)$ for $i=2,3,4$. Let $\cpx{Q_1}$ be the direct
summand $P_A(1)\rightarrow P_A(2)$ of $\cpx{Q}$. Applying $F$ to the
following distinguished triangle in $\Db{A}$
$$P_A(2)[-1]\longrightarrow \cpx{Q_1}[-1]\longrightarrow P_A(1)\longrightarrow P_A(2),$$
we see that $F(P_A(1))$ is of the following form
$$0\longrightarrow P_B(2)\longrightarrow P_B(1)\longrightarrow 0, $$
where $P_B(2)$ is in degree zero. Thus $F(A)$ is isomorphic in
$\Db{B}$ to a complex $\cpx{\bar{Q}}$ which is the direct sum of the
following two complexes
$$\begin{array}{ccccccc} 0 & \lra & P_B(2)                      & \lra & P_B(1)& \lra  & 0,\\
0 & \lra & P_B(2)\oplus P_B(3)\oplus P_B(4)& \lra & 0   &  \lra & 0.
\end{array}$$ Let $G$ be a quasi-inverse of $F$. Then
$\cpx{\bar{Q}}$ is a tilting complex associated to $G$. Clearly,
$_AQ=Q^{-1}=P_A(1)$ and $_B\bar{Q}=\bar{Q}^1=P_B(1)$. It is easy to
see that $F$ satisfies the condition
$\add(_B\bar{Q})=\add(\nu_B\bar{Q})$, but not the condition
$\add(_AQ)=\add(\nu_AQ)$.  Note that $B$ is a Nakayama algebra and
has 16 non-projective indecomposable modules, while $A$ has more
than 16 non-projective indecomposable modules. Thus $A$ and $B$
cannot be stably equivalent.

This example also shows that Corollary \ref{3.10} may be false for
derived equivalences in general. In fact, we have $_AE = P_A(1)$ and
$_B\bar{E}=0$ in this example.

\medskip
Finally, we mention the following questions.

\smallskip
(1) Find new conditions for a derived equivalence to induce a stable
equivalence of Morita type.

\smallskip
(2) Does Theorem \ref{derivedstbaleofMoritatype} hold true for Artin
$R$-algebras $A$ and $B$ such that $A$ and $B$ both are projective
over $R$ ? (For the definition of stable equivalence of Morita type
between Artin algebras see \cite{Xi2}).

\smallskip
(3) Let $F:\Db{A}\longrightarrow\Db{B}$ be a derived equivalence
between Artin algebras $A$ and $B$. If $\add(_AQ)=\add(\nu_AQ)$, is
it true that $\repdim(A)\leq \repdim(B)$ ?

\bigskip
{\bf Acknowledgements.} The research work of C. C. Xi is partially
supported by NSFC (No.10731070).

\medskip
June 12, 2008. Revised: April 15, 2009.

\medskip
Current address of W. Hu: School of Mathematical Sciences, Peking
University, Beijing 100871, P.R.China; email: huwei@math.pku.edu.cn.
\end{document}